\documentclass[12pt]{article}
\usepackage{latexsym,amsfonts,amssymb,amsmath,theorem}

\let\eps=\varepsilon \let\kappa=\varkappa \let\emptyset=\varnothing
\let\cal=\mathcal \let\bbb=\mathbb \let\goth=\mathfrak  
\newcommand\make[2]{\begin{#1}\label{#2}}        
\newcommand\abs[1]{\left|#1\right|}              
\newcommand\norm[1]{\left\|#1\right\|}           
\newcommand\qed{\ifhmode\unskip\nobreak\fi\quad  
   \ifmmode\square\else\hbox{$\square$}\fi}      
\newcommand\proofskip{\vspace{
   \theorempostskipamount}}                      
\input umsa.fd \input umsb.fd \input ueuf.fd     

\theorempreskipamount=\smallskipamount
\theorempostskipamount=\smallskipamount
\newcommand\proof[1]{{\itshape Proof#1\/\,}}
\newcommand\barX{\mspace{4mu}\overline{%
  \mspace{-4mu}X\mspace{-2mu}}\mspace{2mu}}

\newtheorem{theorem}{Theorem}[section]
\newtheorem{proposition}[theorem]{Proposition}

{\theorembodyfont{\upshape}}

\numberwithin{equation}{section}

\makeatletter
\renewcommand*\l@section{\smallskip\@dottedtocline{1}{0em}{1.8em}}
\def\addcontentsline#1#2#3{\addtocontents{#1}{\protect\contentsline{#2}
  {\bfseries#3}{\bfseries\thepage}}}
\makeatother

\begin{document}
\title{\bfseries POSITIVE PROCESSES}
\author{\itshape V.\,I.\,Bakhtin}
\date{}
\maketitle

\vspace{-3\medskipamount}
\centerline{Belarus State University (bakhtin@bsu.by)}

\begin{abstract}
In the present paper we introduce positive flows and processes,
which generalize the ordinary dynamical systems and stochastic
processes. We develop a branch of theory of positive operators
based on the concepts of phase and positive algebras, the
spectral potential, the dual entropy, equilibrium measures, the
action functional, sensitive states, empirical measures and prove
within it the law of large numbers with respect to the sensitive
states and calculate asymptotics for probabilities of large
deviations in terms of the action functional.
\end{abstract}

\bigbreak\medskip
\quad\parbox{14cm}{

\medbreak
}

\bigbreak
\tableofcontents
\bigbreak

\section*{Declaration}  
\addcontentsline{toc}{section}{Declaration}

In this paper we introduce notions of a positive flow and a
positive process. They generalize the ordinary dynamical systems
and stochastic processes. Their main distinctive feature is the
presence of not unique probability measure on the phase space but
a family of measures that is defined as follows. For each
positive flow or process we fix a certain convex functional on
the algebra $C(X)$, where $X$ is the phase space of the object
being studied. This functional is called the spectral potential
and it is none other than the logarithm of the spectral radius of
the corresponding weighted shift operator. It is analogous to the
free energy (the Helmholtz potential) in thermodynamics. Any its
subgradient automatically turns out to be an invariant
probability measure on $X$. Below, as they do in thermodynamics,
we call these measures equilibrium. All of them have equal rights
because there is no reason to give a preference to any of them
rather than to other.

Another peculiarity of the positive flows and processes is the
fact that we treat the concept of motion mostly from the
algebraic point of view. We declare an abstract semigroup of
positive operators (for instance, the family of weighted shift
operators) as the source of any motion. And first of all we
investigate the evolution of operators instead of individual
trajectories in the phase space. This approach is borrowed from
quantum mechanics. Within it the very concept of a phase space
becomes secondary. This tells, in particular, in the fact that
for the same law of motion (evolution) one can consider different
phase spaces. They are determined by the set of those functions
whose values an observer wishes to measure out. In quantum
mechanics similar sets of functions are called Abelian algebras
of observables. But in the present paper we call them the phase
algebras.

The concepts of a positive flow and a positive process are
equivalent in the sense that they are particular cases of each
other. To define them one should specify some semiordered
operator algebra with one distinguished element (evolution
operator) and fix in it some Abelian subalgebra isomorphic to the
algebra of continuous functions on a Hausdorff compact set. The
last set is called the phase space of the positive flow (process).

Despite the great generality of the described approach, we
succeed to prove within it some quite concrete
probability-theoretic results. For instance, if the spectral
potential is differentiable (in the sense of G\^ateaux) at some
point then for the corresponding equilibrium measure the law of
large numbers is valid. Moreover, we obtain exponential
asymptotics for probabilities of observation of various empirical
measures.

Every discrete time dynamical system can be represented as a
positive flow with a fixed equilibrium measure. On the other
hand, if we consider a dynamical system without distinguished
measure on the phase space then we can assign to it many
different positive flows with different families of equilibrium
measures. In the same way positive processes correlate with
stochastic processes.

Finally, notice that positive flows and processes allow us to
introduce in a natural way the ideology and the categories of the
classical equilibrium thermodynamics into the theory of dynamical
systems and stochastic processes (though this is not under
discussion in the text).

The paper is organized as follows. In Section \ref{1..} we
introduce the notion of a dynamical potential. We prove that both
the topological pressure, which is well known in the theory of
dynamical systems, and the logarithm of the spectral radius of a
weighted shift operator are, in fact, dynamical potentials. For
any dynamical potential by means of the Fenchel--Legendre
transform we define the corresponding dual entropy. In Section
\ref{2..} we introduce phase algebras, which are isomorphic to
algebras of continuous functions on the phase space, and
formulate a criterion for an algebra to be phase (a certain analog
of the Gelfand--Naimark theorem). The next Section \ref{3..} has
rather technical character and is entirely devoted to the proof
of the criterion. In Section \ref{4..} we define the base element
of all subsequent constructions --- the positive algebras. They
inherit most important features of operator algebras on the
real spaces $C(X)$ and $L^p(X,\mu)$. It is proved that for any
positive element of a positive algebra there exists a linear
eigenfunctional on the algebra with eigenvalue equal to the
spectral radius of the element. In Section \ref{5..} definitions
of a positive flow and a spectral potential are presented. It is
proved there that every spectral potential is a dynamical
potential. From this fact it follows that the spectral potential
is convex, and that any its subgradient is an invariant
probability measure on the phase space. In Section~\ref{6..} we
define eigen-states of a positive flow and prove the law of large
numbers. In Section \ref{7..} we calculate asymptotics for
probabilities of large deviations by means of an action
functional. Section \ref{8..} deals with positive processes. It
contains a construction that enables us to identify naturally
each positive process with a certain positive flow called the
suspension over the process. On the other hand, by definition
each positive flow is a positive process. So it possesses the
corresponding suspension as well. In Section \ref{9..} we prove
that in essence a positive flow and the suspension over it
coincide. The concluding Section \ref{10..} contains several
examples illustrating the fact that all dynamical systems and
positive processes are special cases of positive flows and
positive processes.

\section{Dynamical potentials and equilibrium measures}\label{1..}

Let $X$ be a Hausdorff compact set. Denote by $C(X)$ the space of
all real-valued continuous functions on~$X$ supplied with the
supremum norm. By the Riesz theorem \cite[Theorem IV.6.3]{Riesz}
each linear continuous functional $\mu\!:C(X)\to\mathbb R$ is
identified with a real-valued Borel measure on $X$, which we will
denote by the same letter $\mu$, so that $\mu(f) =\int_X f\,d\mu$
for all $f\in C(X)$. In particular, if the functional $\mu$ is
positive (i.\,e., nonnegative over nonnegative functions) then
the corresponding measure $\mu$ is positive too. And if the
functional $\mu$ is not only positive but also normalized (takes
the unit value on the unit function) then the corresponding
measure is probability.

Consider an arbitrary continuous mapping $\alpha\!:X\to X$.
Denote by $(X,\alpha)$ the discrete time dynamical system
generated by $\alpha$. For any $\varphi\in C(X)$ define the
weighted shift operator $A_\varphi$ in $C(X)$ by the rule
$A_\varphi f =e^\varphi\cdot f\circ\alpha$. Let
$\lambda(\varphi)$ be the logarithm of the spectral radius of
$A_\varphi$. It is easily shown (which will be done below)
that the functional $\lambda(\varphi)$ possesses a number of
remarkable properties. Firstly, it is {\em monotone:\/} if
$\varphi\le \psi$ then $\lambda(\varphi)\le \lambda(\psi)$.
Secondly, it is {\em additive homogeneous:\/} $\lambda(\varphi+t)
=\lambda(\varphi) +t$ for any $t\in\mathbb R$. Thirdly, $\lambda
(\varphi+\psi) =\lambda(\varphi+\psi\circ\alpha)$ for all
$\varphi,\psi \in C(X)$. This property we shall call the {\em
strong $\alpha$-invariance\/} of $\lambda(\varphi)$. And
fourthly, $\lambda(\varphi)$ is {\em convex\/} with respect to
$\varphi$.

If we fix a positive Borel measure $\mu$ on $X$ then we may
consider the action of the weighted shift operator $A_\varphi$ in
the space $L^p(X,\mu)$, where $p\ge 1$. In this case the four
properties of $\lambda(\varphi)$ mentioned above will be true as
well. And it is easy to check that one more invariant of the
dynamical system $(X,\alpha)$ named the topological pressure
possesses the same properties. So we arrive to the following
definition.

We shall say that a functional $\lambda(\varphi)$ on the $C(X)$
is a {\em dynamical potential\/} of the dynamical system
$(X,\alpha)$ if $\lambda(\varphi)$ is monotone, additive
homogeneous, strongly $\alpha$-invariant, and convex with respect
to $\varphi$. It follows from the definition that any dynamical
potential satisfies the Lipschitz condition with unit Lipschitz
constant. Indeed, the monotonicity and the additive homogeneity
imply that $\lambda(\varphi) -\lambda(\psi) \le \lambda(\psi+\max
\abs{\varphi-\psi}) -\lambda(\psi) = \max\abs{\varphi-\psi}$.

Let us define the {\em dual entropy\/} of a dynamical potential
$\lambda (\varphi)$ as the function of measure
\make{equation}{1,,1}
S(\mu) =\inf_{\varphi\in C(X)}\bigl(\lambda(\varphi) -
\mu(\varphi)\bigr),\qquad \mu\in C^*(X).
\end{equation}
Obviously, $S(\mu)$ is bounded above (because $S(\mu) \le
\lambda(0)$) and only in sign differs from the Fenchel--Le\-gendre
transform of the function $\lambda(\varphi)$. Hence the dual
entropy is always concave and upper semicontinuous (relative to
the *-weak topology). Since the Legendre transform is involution,
we have
\make{equation}{1,,2}
\lambda(\varphi) =\sup_{\mu\in C^*(X)}\bigl(\mu(\varphi) +
S(\mu)\bigr).
\end{equation}

Suppose $\lambda(\varphi)$ is a dynamical potential of a
dynamical system $(X,\alpha)$, where $X$ is compact and
$\alpha\!: X\to X$ is continuous. Any subgradient of the $\lambda
(\varphi)$ at a point $\varphi\in C(X)$ we shall call an {\em
equilibrium measure\/} corresponding to the $\varphi$ (in other
words, a linear functional $\mu\!:C(X)\to\mathbb R$ is an
equilibrium measure iff $\lambda(\varphi +\psi)-
\lambda(\varphi)\ge \mu(\psi)$ for all $\psi\in C(X)$).
Evidently, the set of all equilibrium measures corresponding to
any fixed $\varphi$ is convex and closed (in the *-weak topology).
It follows from the theorem about supporting hyperplane that this
set is nonempty. It consists of only one measure $\mu$ if and
only if there exists the G\^ateaux derivative
$\lambda'(\varphi)$. In this case $\mu =\lambda'(\varphi)$.

From the definitions of the dual entropy and an equilibrium
measure it follows immediately

\make{proposition}{1..1}
For any\/ $\varphi\in C(X)$ and\/ $\mu\in C^*(X)$ the Young
inequality\/ $S(\mu)\le \lambda(\varphi) -\mu(\varphi)$ holds. It
turns into equality if and only if\/ $\mu$ is an equilibrium
measure corresponding to\/ $\varphi$.
\end{proposition}

\make{proposition}{1..2}
If\/ $S(\mu)>-\infty$ then the measure\/ $\mu$ is probability\/
and $\alpha$-invariant. In particular, this is true for all
equilibrium measures.
\end{proposition}

\proof{}
is by contradiction. Firstly, in view of the additive homogeneity
of the dynamical potential, $S(\mu)\le \lambda(t) -\mu(t)
=\lambda(0)+t(1-\mu(1))$ for all $t\in \mathbb R$. And if
$\mu(1)\ne 1$ then $S(\mu) =-\infty$. Secondly, assume that
$\mu(\varphi)<0$ for a certain nonnegative function $\varphi\in
C(X)$. Then, since the dynamical potential is monotone, we see
that $S(\mu)\le\lambda (-t\varphi) -\mu(-t\varphi)\le \lambda(0)+
t\mu(\varphi)$ for all $t>0$ and so $S(\mu) =-\infty$. Thirdly,
the $\alpha$-invariance of $\mu$ is equivalent to the identity
$\mu(\varphi) =\mu(\varphi\circ\alpha)$,\ \ $\varphi\in C(X)$.
Assume that $\mu(\varphi)\ne\mu(\varphi\circ\alpha)$ for a
certain $\varphi\in C(X)$. By the strong $\alpha$-invariance of
the dynamical potential, we have
$\lambda(t\varphi -t\varphi\circ\alpha) \equiv \lambda(0)$.
Therefore, $S(\mu)\le\lambda(t\varphi-t\varphi \circ \alpha)
-\mu(t\varphi -t\varphi\circ\alpha) =\lambda(0)-t (\mu(\varphi) -
\mu(\varphi \circ\alpha))$ and hence $S(\mu) =-\infty$.
Finally, if $\mu$ is an equilibrium measure corresponding to
$\varphi\in C(X)$ then by Proposition \ref{1..1} we have
$S(\mu) =\lambda(\varphi) -\mu(\varphi) >-\infty$.
\qed

\begin{proposition}\label{1..3}
If\/ $S(\mu) >-\infty$ then\/ $\mu$ belongs to the closure of all
equilibrium measures\/ \textup{(}in the norm of\/
$C^*(X)$\textup{)}.
\end{proposition}

\proof.
A variant of Bishop--Phelps' theorem \cite{Bishop} states that if
$\lambda(\varphi)$ is a continuous convex functional on a Banach
space $B$ and $\mu\in B^*$ and the difference $\lambda- \mu$ is
bounded below then $\mu$ lies at zero distance from the set of
subgradients of $\lambda$.
\qed
\proofskip

Let us define on $C^*(X)$ the {\em action functional\/}
$\,\tau_\varphi(\mu) =\mu(\varphi)+S(\mu)- \lambda(\varphi)$. By
Proposition~\ref{1..1} it is always nonpositive and it vanishes
only on the equilibrium measures. And, in view of
Proposition~\ref{1..3}, it may be finite (not equal to $-\infty$)
only on the closure of all equilibrium measures.

Denote by $M_\alpha(X)$ the set of all $\alpha$-invariant
Borel probability measures on $X$. From Proposition~\ref{1..2} it
follows that the supremum in \eqref{1,,2} is attained on
$M_\alpha(X)$. Hence any dynamical potential has the form
\make{equation}{1,,3}
\lambda(\varphi) =\sup_{\mu\in M_\alpha(X)}
\bigl(\mu(\varphi) +S(\mu)\bigr),
\end{equation}
where $S(\mu)$ is the corresponding dual entropy. By Proposition
\ref{1..1}, the supremum in \eqref{1,,3} is attained exactly on
the equilibrium measures corresponding to the $\varphi$.
Therefore the uniqueness of the extremal measure in \eqref{1,,3}
is equivalent to the existence of the G\^ateaux derivative
$\lambda'(\varphi)$.

Now we present one more statement explaining the structure of
dynamical potentials.

\make{proposition}{1..4}
If a function\/ $S(\mu)$ is bounded above on\/ $M_\alpha(X)$
then\/ \eqref{1,,3} determines a dynamical potential.
\end{proposition}

Proof of this fact reduces to the trivial verification of all
conditions defining a dynamical potential and so we leave it to
the reader.

\proofskip
Since dual entropy satisfies $S(\mu)\le \lambda(0)$, we see
from \eqref{1,,3} that $\lambda(\varphi) \le
\sup_{M_\alpha(X)}\mu(\varphi) +\lambda(0)$. On the other hand,
if $S(\mu) \equiv \lambda(0)$ then by means of \eqref{1,,3} we
obtain the dynamical potential of the form $\lambda(\varphi)
=\sup_{M_\alpha(X)} \mu(\varphi) +\lambda(0)$. The previous
estimate shows that this $\lambda(\varphi)$ is maximal among
all dynamical potentials with fixed value of $\lambda(0)$.

\section{Phase algebras}\label{2..}

Let $L$ be a vector space over $\mathbb R$. Recall that a {\em
cone\/} in $L$ is a subset $K\subset L$ that is convex,
homogeneous (if $v\in K$ and $t\ge 0$ then $tv\in K$), and
antisymmetric (if $v\in K$ and $-v\in K$ then $v=0$). Every cone
$K$ determines a partial order on $L$: we say that $u\le v$ iff
$v-u\in K$. In particular, vectors from the cone are called {\em
positive.} A real vector space with fixed cone is said to be
{\em semiordered.} A linear functional on the semiordered space is
called positive if it is nonnegative on the cone. A linear
operator on the semiordered space is called positive if it maps
the cone into itself.

A {\em semiordered algebra\/} is an algebra $\mathcal A$ with
fixed cone $\mathcal A_+ \subset \mathcal A$ such that the
product of any two positive elements is also positive. For
example, if $\mathcal A$ is an arbitrary algebra of functions on
a set $X$ then it is naturally semiordered by the cone of
nonnegative functions. In particular, for any compact set $X$ the
corresponding algebra $C(X)$ is semiordered. Evidently, each
subalgebra $\mathcal B$ of a semiordered algebra $\mathcal A$ is
semiordered too (by the cone $\mathcal B\cap\mathcal A_+$).

A semiordered algebra $\mathcal C$ with unit we shall call a {\em
phase algebra\/} if there exists a Hausdorff compact set~$X$ such
that $\mathcal C$ is isomorphic (including the order) to some
dense (relative to the supremum norm) subalgebra in~$C(X)$. This
$X$ we shall call the {\em support\/} or the {\em phase space\/}
of the phase algebra $\mathcal C$.

\make{proposition}{2..1}
The support of a phase algebra is defined uniquely up to a
homeomorphism.
\end{proposition}

\proof.
Assume that a phase algebra $\mathcal C$ is isomorphic (including
the order) to some dense subalgebras $\cal C_X\subset C(X)$ and
$\cal C_Y\subset C(Y)$, where $X$ and $Y$ are compact. Then
$\mathcal C_X$ is isomorphic to $\mathcal C_Y$ including the
order. Under this isomorphism the unit function has to go to the
unit function. Hence every constant function goes to the same
constant function. But since the isomorphism between $\mathcal
C_X$ and~$\mathcal C_Y$ preserves the order, it follows that it
preserves the supremum norm. By continuity it extends up to the
metric isomorphism between $C(X)$ and $C(Y)$. By the
Gelfand--Naimark theorem each compact set $X$ is homeomorphic to
the space of maximal ideals in $C(X)$. Therefore the isomorphism
bertween $C(X)$ and $C(Y)$ generates a homeomorphism between $X$
and $Y$. \qed

\make{proposition}{2..2}
Every positive linear functional on a phase algebra is generated
by a certain Borel measure on the support of the algebra.
\end{proposition}

\proof.
Without loss of generality it is sufficient to consider a
positive linear functional $\mu$ on a phase algebra $\mathcal C$
dense in $C(X)$. In this case for any $f\in\mathcal C$ we have
$\,\mu(f)\le \mu(\max f) =\mu(1)\max f\,$ and, on the other hand,
$\,\mu(f)\ge \mu(\min f) =\mu(1)\min f$. So $\mu$ is continuous
relative to the supremum norm on $\mathcal C$. By the continuity
it extends to $C(X)$ preserving the positiveness. But by the
Riesz theorem each positive linear functional $\mu$ on $C(X)$ is
identified with a certain Borel measure $\mu$ on $X$ so that
$\mu(f)\equiv \int_X f\,d\mu$.
\qed

\make{proposition}{2..3}
Suppose\/ $\mu$ is a positive linear functional on a phase
algebra\/ $\mathcal C$. Then the functional\/ $\lambda(\varphi) =
\ln\mu(e^\varphi)$ is convex with respect to\/
$\varphi\in{\mathcal C}$.
\end{proposition}

\proof.
Without loss of generality it can be assumed that $\mathcal C$ is
a dense subalgebra in $C(X)$ and $\mu$ is a Borel measure on $X$.
So the expression $\mu(e^\varphi)$ is well defined even if
$e^\varphi\notin\cal C$. It is easy to see that the function
$\mu(fg)$ of the pair $(f,g)\in{\cal C}^2$ is bilinear, symmetric,
and positive definite. Therefore it satisfies the
Cauchy--Bunyakovskii inequality $\mu(fg)^2\le \mu(f^2)\mu(g^2)$.
To prove the convexity of $\lambda(\varphi)$ it is enough to
check that the restriction of $\lambda(\varphi)$ to any straight
line $\{\varphi_t = \varphi_0 +t\varphi\}$ is convex with respect
to~$t\in\bbb R$. Let us differentiate this restriction twice in
$t$:
\begin{equation*}
\frac{d\lambda(\varphi_t)}{dt} =
\frac{\mu(e^{\varphi_t}\varphi)}{\mu(e^{\varphi_t})},\qquad
\frac{d^2\lambda(\varphi_t)}{dt^2} =
\frac{\mu(e^{\varphi_t}\varphi^2)\mu(e^{\varphi_t}) -
\mu(e^{\varphi_t}\varphi)^2}{\mu(e^{\varphi_t})^2}.
\end{equation*}
By the Cauchy--Bunyakovskii inequality, $\mu\left(e^{\varphi_t/2}
\cdot e^{\varphi_t/2}\varphi\right)^2\le \mu(e^{\varphi_t})\mu
(e^{\varphi_t}\varphi^2)$. It is seen from here that $d^2\lambda
(\varphi_t) /dt^2\ge 0$. Hence $\lambda(\varphi)$ is convex.
\qed
\proofskip

Now consider the question what kind of algebras are phase. One
important class of phase algebras is provided by

\make{proposition}{2..4}
Suppose\/ $Y$ is a Hausdorff compact set and\/ $\cal I$ is a
closed ideal in\/ $C(Y)$. Then the quotient algebra\/ $C(Y)/\cal
I$ is phase.
\end{proposition}

\proof.
Denote by $X$ the annihilator of $\cal I$. In other words, $X
=\bigcap_{f\in\cal I}f^{-1}(0)$. This $X$ is closed and so
compact. The restriction of functions from $Y$ to $X$ generates
the algebraic homomorphism $\pi\!:C(Y)\to C(X)$. By the
Brauer--Urysohn theorem any function $f\in C(X)$ can be extended
up to a continuous function on $Y$. Therefore $\pi$ is a
surjection. On the other hand, any closed ideal $\cal I \subset
C(Y)$ coincides with the set of functions vanishing on its
annihilator $X$. Hence $\cal I$ is the kernel of $\pi$. From here
it follows that $C(Y)/\cal I \cong C(X)$.
\qed
\proofskip

In fact we are able to formulate a necessary and sufficient
condition for a semiordered algebra to be phase. To this end
consider a semiordered algebra $\cal C$ with cone $\cal C_+$ and
positive unit $I\in\cal C_+$. An element $a\in\cal C$ we shall
call {\em bounded\/} if there exists a number $t>0$ such that
$-tI\le a\le tI$. From the definition of the phase algebra we see
that it consists of bounded elements and is commutative
and the intersection of the positive cone with any affine
straight line is closed. It turns out that the inverse is valid
as well.

\make{theorem}{2..5}
If a semiordered Abelian algebra\/ $\cal C$ has a
positive unit\/ $I\in\cal C_+$ and consists of bounded elemtnts
and the intersection of\/ $\cal C_+$ with any affine straight
line is closed, then the algebra\/ $\cal C$ is phase.
\end{theorem}

This theorem can be treated as a real analog of
Gelfand--Naimark's theorem about structure of an Abelian
$C^*$-algebra with unit \cite{Gelfand}. It will be proved in the
next section. Theorem \ref{2..5} provides at least two important
examples of phase algebras. First of them is any algebra $\cal C$
consisting of bounded real-valued functions on an arbitrary
set~$X$ and containing the identical unit. The second is the set
of essentially bounded functions $L^\infty(X,\mu)$ on any measure
space $(X,\mu)$.

\section[Proof of the real analog of Gelfand--Naimark's theorem]
{Proof of the real analog of Gelfand--Naimark's\\
theorem}\label{3..}

Consider the algebra $\cal C$ from Theorem \ref{2..5}. Set
\make{equation}{3,,1}
\norm a =\inf\{\,t>0\mid -tI\le a\le tI\,\},\qquad a\in\cal C.
\end{equation}

\make{lemma}{3..1}
Under the conditions of Theorem \ref{2..5} the function\/
\eqref{3,,1} is a norm on\/ $\cal C$. The completion of the
cone\/ $\cal C_+$ in this norm is a cone as well (in the
completion of\/ $\cal C$).
\end{lemma}

\proof.
If we regard $\cal C$ as a vector space then, evidently,
\eqref{3,,1} is a seminorm on it. Assume that $\norm a=0$. In
this case it follows from \eqref{3,,1} that $tI+a\in\cal C_+$ and
$tI-a\in\cal C_+$ for all $t>0$. Then by the condition of Theorem
\ref{2..5} we have $a\in\cal C_+$ and $-a\in\cal C_+$. Therefore
$a=0$. This proves that $\norm a$ is a norm on the vector space
$\cal C$.

Further, assume that $\norm a <t$ and $\norm b <s$. Then $tI\pm
a\in \cal C_+$ and $sI\pm b\in\cal C_+$. Therefore,
\begin{align*}
2tsI -2ab &= (tI-a)(sI+b) +(tI+a)(sI-b) \in\cal C_+, \\
2tsI +2ab &= (tI+a)(sI+b) +(tI-a)(sI-b) \in\cal C_+.
\end{align*}
So we see that $-tsI\le ab\le tsI$. Hence $\norm{ab} \le\norm
a\norm b$ for all $a,b\in\cal C$. This proves that
$\norm{\,\cdot\,}$ is a norm on the algebra $\cal C$.

Let $\left[\cal C_+\right]$ be the completion of $\cal C_+$ in
the norm. Obviously, $\left[\cal C_+\right]$ is convex and
homogeneous. By \eqref{3,,1}, we have $\norm a\le\norm{a+b}$ for
all $a,b\in \cal C_+$. This inequality still remains in power for
$a,b \in \left[\cal C_+\right]$. Therefore, if $a\in \left[\cal
C_+\right]$ and at the same time $-a\in \left[\cal C_+\right]$
then $\norm a\le\norm{a-a} =0$. Hence $\left[\cal C_+ \right]$ is
a cone.
\qed

\make{lemma}{3..2}
For each natural\/ $n$,
\begin{equation*}
x^2 =4^{-n}\sum_{i=-n}^n \frac{i^2}{n^2}C_{2n}^{n+i} (1+x)^{n+i}
(1-x)^{n-i} -\frac{1-x^2}{2n}.
\end{equation*}
\end{lemma}

\proof{}
is by direct calculation. Note that
\begin{equation*}
4^n =\bigl((1+x)+(1-x)\bigr)^{2n} =\sum_{i=-n}^n C_{2n}^{n+i}
(1+x)^{n+i}(1-x)^{n-i}.
\end{equation*}
Therefore,
\begin{align*}
4^n -\sum_{i=-n}^n\frac{i^2}{n^2} &C_{2n}^{n+i} (1+x)^{n+i}(1-x)^{n-i}
=\sum_{i=-n}^n\left(1-\frac{i^2}{n^2}\right) C_{2n}^{n+i} (1+x)^{n+i}
(1-x)^{n-i} =\\[3pt]
&= \sum_{i=-n+1}^{n-1} \frac{(n+i)(n-i)}{n^2}\frac{(2n)!}{(n+i)!
(n-i)!} (1+x)^{n+i}(1-x)^{n-i} =\\[3pt]
&= \frac{2n(2n-1)}{n^2}(1-x^2)\sum_{i=-n+1}^{n-1} C_{2n-2}^{n+i-1}
(1+x)^{n+i-1}(1-x)^{n-i-1} =\\[3pt]
&= \left(4-\frac{2}{n}\right)(1-x^2)\,4^{n-1} =4^n\left(1-x^2-
\frac{1-x^2}{2n}\right). \qed
\end{align*}
\vspace{3pt}

Denote by $[\cal C]$ the completion of $\cal C$ in the norm
\eqref{3,,1} and by $\left[\cal C_+\right]$ the closure of $\cal
C_+$ in the same norm. For each $a\in\cal C$ define two numbers:
$\sup a=\min\{\,t\in\bbb R\mid a\le tI\,\}$ and $\inf a =\max\{\,
t\in\bbb R\mid a\ge tI\,\}$, and let $I_a= [\inf a,\sup a]$.

\make{lemma}{3..3}
Under the conditions of Theorem \ref{2..5} for each\/ $a\in\cal
C$ there exists a unique algebraic homomorphism\/ $\pi_a\!:
C(I_a)\to [\cal C]$ that maps nonnegative functions to\/ $\left[
\cal C_+\right]$ and takes any polynomial\/ $P(x) =c_0+c_1x
+\dots+c_nx^n$ to\/ $P(a) =c_0I+c_1a+\dots+c_na^n$.
\end{lemma}

\proof.
First let us check that the square of any element $b\in\cal C$
is positive. Without loss of generality, it can be assumed that
$-I\le b\le I$. Then by Lemma \ref{3..2},
\begin{equation*}
b^2+\frac{I-b^2}{2n} =4^{-n}\sum_{i=-n}^n\frac{i^2}{n^2}
C_{2n}^{n+i}(I+b)^{n+i}(I-b)^{n-i}\in \cal C_+.
\end{equation*}
By condition, the intersection of $\cal C_+$ with any affine
straight line is closed. Therefore, we get $b^2\in\cal C_+$ as
$n\to\infty$.

Thereupon let us prove that if a polynomial $P(x)$ is nonnegative
on the segment $I_a =[\inf a,\sup a]$, then $P(a) \in\cal C_+$.
Indeed, in this case there exists a factorization
$P(x) =c\prod_i(x-\alpha_i)\*\prod_j
(\beta_j-x)\*\prod_k\bigl((x-\gamma_k)^2+\delta_k^2\bigr)$ in which
$c\ge 0$,\ \ $\alpha_i \le \inf a$, and $\beta_j\ge \sup a$.
Hence, $P(a) =c\prod_i (a-\alpha_iI)\*\prod_j
(\beta_jI-a)\*\prod_k\bigl((a- \gamma_kI)^2+ \delta_k^2I\bigr)$ is
a product of factors lying in $\cal C_+$ and so it also belongs
to $\cal C_+$.

From here it follows that if $-t\le P(x)\le t$ on the segment
$I_a$ for some $t>0$ then $-tI\le P(a)\le tI$. Thus
$\norm{P(a)}\le \sup_{I_a} \abs{P(x)}$. This proves continuity of
the mapping $P(x)\mapsto P(a)$. By the continuity this mapping
extends up to a positive algebraic homomorphism
$\pi_a\!:C(I_a)\to [\cal C]$.
\qed
\proofskip

{\em Proof of Theorem \ref{2..5}.\,}
For each $a\in\cal C$ define the segment $I_a =[\inf
a,\allowbreak\sup a]$ and the corresponding homomorphism
$\pi_a\!: C(I_a)\to [\cal C]$ from Lemma \ref{3..3}. Consider the
Cartesian product $Y =\prod_{a\in \cal C}I_a$. By Tychonoff's
theorem it is compact in the cylinder topology. Let $\cal A$ be
the subalgebra in $C(Y)$ consisting of various finite sums of
products $\sum_{i=1}^m f^i_1(x_1)\dotsm f^i_n(x_n)$, where
$x_j\in I_{a_j}$ and $f^i_j\in C(I_{a_j})$. This subalgebra
contains the identical unity and separates points of $Y$. By
Stone's theorem $\cal A$ is dense in $C(Y)$. Define an algebraic
homomorphism $\pi\!:\cal A\to[\cal C]$ by means of the formula
$\pi\left(\sum_{i=1}^m f^i_1(x_1)\dotsm f^i_n(x_n)\right)
=\sum_{i=1}^m \pi_{a_1}(f^i_1)\dotsm\pi_{a_n}(f^i_n)$, where
$\pi_{a_j}$ are the homomorphisms from Lemma \ref{3..3}.

Assume that a function $F\in\cal A$ is nonnegative. Then for any
$\eps>0$ the sum $F+\eps$ may be represented as $F+\eps
=\sum_{i=1}^k g^i_1\dotsm g^i_n$, where all $g^i_j\in C(I_{a_j})$
are nonnegative. This will be proved below in Theorem~\ref{3..4}.
Combining this representation with Lemma \ref{3..3}, we obtain
$\pi(F+\eps) \in\left[\cal C_+\right]$, which implies
$\pi(F)\in\left[\cal C_+\right]$ as $\eps\to 0$. Hence the
homomorphism $\pi$ is positive. If a function $F\in\cal A$
satisfies $-t\le F\le t$,\ \ $t>0$, then $-tI\le \pi(F)\le tI$.
Therefore $\norm{\pi(F)}$ does not exceed the supremum norm of
the $F$. By continuity $\pi$ uniquely extends from $\cal A$ to
$C(Y)$ preserving positiveness.

Denote by $\cal I$ the kernel of the homomorphism
$\pi\!:C(Y)\to[\cal C]$. Obviously, $\cal I$ is a closed ideal
in~$C(Y)$. Let $X$ be its annihilator: $X=\bigcap_{f\in\cal I}
f^{-1}(0)$. The set $X$ is closed and so compact. Define a
homomorphism $\rho\!: C(X)\to [\cal C]$ in the following way. By
the Brauer--Urysohn theorem any $f\in C(X)$ can be extended up to
a continuous function $F$ on $Y$. Set $\rho(f) = \pi(F)$. It is
well known that every closed ideal $\cal I\subset C(Y)$ coincides
with the set of functions vanishing on its annihilator $X$.
Therefore the homomorphism $\rho$ is well defined and injective.

For any $b\in\cal C$ define the function $F_b\in C(Y)$ by the
rule $F_b\left(\{x_a\in I_a\}_{a\in\cal C}\right) =x_b$. Let
$f_b$ be the restriction of $F_b$ to $X$. Then $\rho(f_b)
=\pi(F_b) =b$ and hence $\cal C\subset \rho(C(X))$. On the other
hand, the functions $f_b$, where $b\in\cal C$, separate points of
$X$. So by Stone's theorem the algebra $\rho^{-1}(\cal C)$ is
dense in $C(X)$. Since the homomorphism $\rho$ is injective, its
inverse $\rho^{-1}\!: b \mapsto f_b$ establishes an isomorphism
between $\cal C$ and $\rho^{-1} (\cal C)$. It is easily seen that
this isomorphism is positive. Indeed, by construction the
functions $F_b$ and $f_b$ take values in the segment $I_b =[\inf
b,\sup b]$. And if $b\in\cal C_+$ then $\inf b\ge 0$. This
completes the proof of Theorem \ref{2..5}.
\qed

\make{theorem}{3..4}
Suppose on the Cartesian product\/ $X_1\times\dots\times X_n$ of
Hausdorff compact sets\/ $X_i$ there is a nonnegative function\/
$F$ of the form\/ $F=\sum_{i=1}^m f^i_1\dotsm f^i_n$, where\/
$f^i_j\in C(X_j)$. Then for any\/ $\eps>0$ there exist
nonnegative functions\/ $g^i_j\in C(X_j)$ such that\/ $F+\eps
=\sum_{i=1}^k g^i_1\dotsm g^i_n$.
\end{theorem}

Let us denote by $\cal A_\#$ the set of all finite sums of the form
$\sum_{i=1}^k g^i_1\dotsm g^i_n$, where all $g^i_j\in C(X_j)$ are
nonnegative. To prove Theorem \ref{3..4} we need the following

\make{lemma}{3..5}
If under the conditions of Theorem \ref{3..4} the function\/ $F$
has the form\/ $F=f_1\dotsm f_n$, where\/ $f_j\in C(X_j)$, and\/
$\eps\ge \norm{f_1}\times\dots\times\norm{f_n}$, then\/
$F+\eps\in\cal A_\#$.
\end{lemma}

\proof.
Without loss of generality it can be assumed that $\eps=
\norm{f_1} \times \dots\times\norm{f_n}$. Set $f^+_j
=\max\{0,f_j\}$ and $f^-_j =\max\{0,-f_j\}$. Then $f_j =f^+_j
-f^-_j$. Accordingly, the $F$ is presented as a sum of $2^n$
terms of the form $\pm f_1^\pm\dotsm f_n^\pm$. Obviously, the
functions $f^+_j,\,f^-_j$ and $g_j =\norm{f_j}- f^+_j -f^-_j$ are
nonnegative and continuous. Hence by distributivity
$\eps=(f^+_1+f^-_1+ g_1)\times\dots\times(f^+_n+f^-_n+g_n)$ is a
sum of $3^n$ terms, each of them being a product of $n$
nonnegative functions. And all the products $f^\pm_1\dotsm
f^\pm_n$ are contained among these $3^n$ terms. Therefore
$F+\eps\in\cal A_\#$.
\qed\proofskip

{\em Proof of Theorem \ref{3..4}.\,}
Take any $\eps>0$. For each $j=1,\,\dots,\,n$ choose a continuous
partition of unity $\varphi^1_j$, \dots, $\varphi^{N_j}_j$ on
$X_j$ such that the oscillations of the functions $f^i_j$ on the
supports of $\varphi^k_j$ do not exceed $\eps$. In the support of
each $\varphi^k_j$ select a point $x^k_j\in X_j$. Put
\begin{equation*}
\eta_\eps(x) =
\begin{cases}
-\eps,& x<-\eps\\ x,& -\eps\le x\le \eps\\ \eps,& x>\eps.
\end{cases}
\end{equation*}
Let us consider the functions $g^{ik}_j(x) =\eta_\eps\left(
f^i_j(x) -f^i_j\bigl(x^k_j\bigr)\right)$ on $X_j$. By
construction $\norm{g^{ik}_j} \le\eps$ and $f^i_j(x) =g^{ik}_j(x)
+f^i_j\bigl(x^k_j\bigr)$ on the $\mathop{\mathrm{supp}}
\varphi^k_j$. Therefore $f^i_j\equiv\sum_{k=1}^{N_j}
\left(g^{ik}_j +f^i_j\bigl(x^k_j\bigr)\right)\varphi^k_j$.

Suppose $t=\max_{i,j}\norm{f^i_j}$. Let us fix indices
$k_1,\,\dots,\, k_n$ and consider the product
$\bigl(g^{ik_1}_1+f^i_j\bigl(x^{k_1}_1\bigr)\bigr)\times\cdots
\times\bigl(g^{ik_n}_n +f^i_n \bigl(x^{k_n}_n\bigr)\bigr)$.
Decompose it by distributivity into the sum of $2^n$ terms and
apply Lemma \ref{3..5} to each of them except the last one. As a
result we see that
\begin{equation*}
\prod_{j=1}^n\left(g^{ik_j}_j+f^i_j\bigl(x^{k_j}_j\bigr)\right)
-\prod_{j=1}^n f^i_j\bigl(x^{k_j}_j\bigr) +(2^n-1)\eps
(t+\eps)^{n-1}\in \cal A_\#.
\end{equation*}
Multiply the last expression by $\varphi^{k_1}_1 \dotsm
\varphi^{k_n}_n$ and sum up over all indices $i$, $k_1,\,
\dots,\,k_n$. By construction the sum will lie in $\cal A_\#$:
\begin{gather*}
\sum_{i=1}^m\prod_{j=1}^n\sum_{k_j=1}^{N_j}\left(g^{ik_j}_j +
f^i_j\bigl(x^{k_j}_j\bigr)\right)\varphi^{k_j}_j -
\sum_{i=1}^m\sum_{k_1\dots k_n}\prod_{j=1}^n f^i_j
\bigl(x^{k_j}_j\bigr)\varphi^{k_j}_j\\[3pt]
+(2^n-1)\eps(t+\eps)^{n-1}\sum_{i=1}^m\prod_{j=1}^n
\sum_{k_j=1}^{N_j}\varphi^{k_j}_j \\[3pt]
=\sum_{i=1}^m f^i_1\dotsm f^i_n -\sum_{k_1\dots k_n}\sum_{i=1}^m
f^i_1\bigl(x^{k_1}_1\bigr)\dotsm f^i_n\bigl(x^{k_n}_n\bigr)
\varphi^{k_1}_1\dotsm\varphi^{k_n}_n +(2^n-1)\eps(t+\eps)^{n-1}m
\\[3pt]
=F-\sum_{k_1\dots k_n}F\bigl(x^{k_1}_1,\dots,x^{k_n}_n\bigr)
\varphi^{k_1}_1\dotsm \varphi^{k_n}_n +(2^n-1)\eps(t+\eps)^{n-1}m
\in\cal A_\#.
\end{gather*}
Since $F$ and $\varphi^k_j$ are nonnegative and $\eps>0$ is
arbitrary this implies Theorem \ref{3..4}.
\qed

\section{Positive algebras}\label{4..}

Assume that a semiordered vector space $L$ with cone $K$ is
supplied with a norm. We shall call this norm {\em monotone\/} if
$\norm u\le\norm{u+v}$ for all $u,v\in K$. The semiordered vector
space with monotone norm will be called {\em positive}.

At once note that the definition of the positive space $L$ stays
in power if we replace the positive cone~$K$ by its closure
$[K]$. In this case the closure of $K$ will stay antisymmetric:
if $v\in [K]$ and $-v\in [K]$ then $\norm v\le\norm{v-v} =0$ (by
monotonicity of the norm).

Examples of the positive spaces are the spaces of integrable
functions $L^p(X,\mu)$, where $1\le p\le\infty$, and the spaces
of summable sequences $l^p$. On the other hand, the space of
differentiable functions with the cone of nonnegative functions
is not positive.

\make{proposition}{4..1}
If\/ $v$ is a positive vector in a positive space\/ $L$ then
there exists a positive linear functional\/ $\nu\!: L\to
\mathbb R$ such that\/ $\norm\nu =1$ and $\nu(v) =\norm v$.
\end{proposition}

\proof.
Without loss of generality it can be assumed that $\norm v=1$.
Denote by $K$ the positive cone in~$L$. Consider two convex sets:
the open unit ball $B=\{u\in L: \norm u<1\}$ and the shifted cone
$K_v= v+K$. From the definition of the positive space it follows
that $B\cap K_v =\emptyset$. By the separability theorem there
exists a linear functional $\nu$ on $L$ such that $\nu(u) \le 1$
for all $u\in B$ and $\nu(u)\ge 1$ for all $u\in K_v$. Obviously,
$\norm\nu =\nu(v) =1$. Besides, for any $u\in K$ we have $\nu(u)
=\nu(v+u)-\nu(v)\ge 0$.
\qed\proofskip

Recall that an algebra $\cal A$ is said to be normed if $\cal A$
is a normed linear space and $\norm{ab}\le \norm a\norm b$ for
all $a,b\in\cal A$. If $\cal A$ contains the unity $I$ then in
addition $\norm I =1$ is required (which always can be provided
by mere replacement of the norm by an equivalent one).

\proofskip
\mbox{\bfseries Definition.}
We shall call a real algebra $\cal A$ {\em positive\/} if it is
normed, semiordered, has a positive unit, and the norm on $\cal
A$ is monotone.

\make{example}{ex4..1}
Any algebra of functions on an arbitrary domain that contains the
identical unity and is supplied with the $\sup$ norm is positive.
\end{example}

\make{example}{ex4..2}
Let $X$ be a topological space. Then the algebra of bounded
linear operators on $C(X)$ is positive. Indeed, if two linear
operators $A,\,B$ on $C(X)$ are positive (take nonnegative
functions to nonnegative ones) then $\norm A =\sup A(1)\le
\sup\bigl(A(1) +B(1)\bigr) =\norm{A+B}$.
\end{example}

\make{example}{ex4..3}
Let $(X,\mu)$ be a measure space. Then the algebra of bounded
linear operators on the space $L^p(X,\mu)$, where $1\le p\le
\infty$, is positive.
\end{example}

\proof.
For any $f\in L^p(X,\mu)$ define the functions $f^+
=\max\{f,0\}$ and $f^- =\max\{-f,0\}$. Then $f=f^+ -f^-$ and
$\abs f =f^+ +f^-$. Suppose $A$ is a positive linear operator on
$L^p(X,\mu)$. Obviously, $Af^+\ge 0$ and $Af^-\ge 0$, and hence
$\abs{Af} =\abs{Af^+ -Af^-}\le Af^+ +Af^- = A\abs{f}$. On the
other hand $\norm f =\bigl\|\abs f\bigr\|$. Therefore,
\begin{equation*}
\norm A =\sup_{f\in L^p(X,\mu)}\frac{\norm{Af}}{\norm f} =
\sup_{f\ge 0}\frac{\norm{Af}}{\norm f}.
\end{equation*}
From the last equality we see that if operators $A$ and $B$ are
positive then $\norm A\le\norm{A+B}$.
\qed

\make{example}{ex4..4}
Assume that $L$ is a semiordered and normed vector space and
$\cal A_+$ is the set of all positive bounded linear operators on
$L$. Denote by $\cal A$ the linear hull of $\cal A_+$. In other
words, $\cal A =\{\,A-B \mid A,B \in\cal A_+\,\}$. For each $A\in
\cal A$ define the seminorm
\begin{equation*}
\abs A =\inf\bigl\{\norm{A'}\bigm| A'\pm A\in\cal A_+\bigr\} =
\inf\bigl\{\norm{A'}\bigm| -A'\le A\le A'\bigr\}.
\end{equation*}
Let $\cal I =\bigl\{\, A\in \cal A\bigm| \abs A =0\,\bigr\}$.
Then the quotient algebra $\cal A/\cal I$ with norm
$\abs{\,\cdot\,}$ is positive.
\end{example}

\proof.
If $A,B\in\cal A_+$ then $\abs A\le \abs{A+B}$ by the definition
of the seminorm. Assume that for some $A,B,A',B'\in\cal A$ we
have $A'\pm A\in\cal A_+$ and $B'\pm B\in \cal A_+$. Then
\begin{align*}
2A'B'-2AB &=(A'-A)(B'+B) +(A'+A)(B'-B)\in\cal A_+,\\[1pt]
2A'B'+2AB &=(A'+A)(B'+B) +(A'-A)(B'-B)\in\cal A_+.
\end{align*}
Hence $A'B'\pm AB\in\cal A_+$. It follows from here that
$\abs{AB} \le \abs A\abs B$ for all $A,B\in\cal A$. Therefore
$\abs{\,\cdot\,}$ is a seminorm on the algebra $\cal A$ and is a
norm on the algebra $\cal A/\cal I$.
\qed

\make{example}{ex4..5}
We shall say that a norm on a semiordered vector space $L$ is
pseudomonotone if there exists a positive constant $C$ such that
$\norm v \le C\norm u$, provided $u\ge v\ge 0$. It is known that
for any pseudomonotone norm there exists an equivalent monotone
norm. It can be given by $\norm u =\sup_{\nu\in V}\abs{\nu(u)}$,
where $V$ is the set of positive linear functionals on $L$ with
unit norm (see \cite{Krasn}). Unfortunately, in the case of
semiordered algebras the estimate $\norm{uv}\le \norm u\norm v$
may fail for this norm. So it would be interesting to know
whether there is a monotone algebraic norm equivalent to the
pseudomonotone norm on a semiordered algebra.
\end{example}

Suppose $\cal A$ is a positive algebra with closed cone $\cal
A_+$. Obviously, every subalgebra $\cal C\subset \cal A$
containing the unity is also positive (with the cone $\cal C_+
=\cal C\cap \cal A_+$). A subalgebra $\cal C\subset\cal A$ will
be called below {\em phase\/} if it is commutative, contains the
unity of $\cal A$, and consists of bounded (in the sense of
definition given in Section \ref{2..}) elements. From Theorem
\ref{2..5} it follows that any phase subalgebra of a positive
algebra $\cal A$ with closed cone $\cal A_+$ is phase in the
sense of definition from Section \ref{2..} as well. In other
words, it is isomorphic to a dense subalgebra in $C(X)$, where
$X$ is a Hausdorff compact set, which is uniquely determined by
the phase algebra and named its support.

\make{example}{ex4..6}
Let $\cal A$ be the set of bounded linear operators on $C(Y)$,
where $C(Y)$ is the algebra of bounded continuous functions on a
certain topological space $Y$. Then $C(Y)$ can be considered as a
phase subalgebra of $\cal A$ (under the condition that each $f\in
C(Y)$ is identified with the operator of multiplication by $f$).
\end{example}

\make{example}{ex4..7}
Supplse $\cal A$ is the set of bounded linear operators on
$L^p(Y,\mu)$, where $(Y,\mu)$ is a measure space. Then
$L^\infty(Y,\mu)$ is a phase subalgebra of $\cal A$.
\end{example}

\make{proposition}{4..2}
If\/ $\cal A$ is a complete in norm positive algebra with closed
cone\/ $\cal A_+$ and\/ $\cal C$ is a phase subalgebra in it then
for any\/ $a\in\cal A_+$ the functional\/ $\lambda(\varphi)
=\ln\norm{e^\varphi a}$ is convex with respect to\/
$\varphi\in\cal C$.
\end{proposition}

\proof.
The phase subalgebra $\cal C$ is isomorphic to a dense subalgebra
in $C(X)$, where $X$ is the support of $\cal C$. This isomorphism
extends up to a positive continuous embedding of $C(X)$ into
$\cal A$. So it follows that for any $\varphi\in\cal C$ the
element $e^\varphi\in \cal A$ is positive. Let $V$ be the set of
all positive linear functionals on $\cal A$ with unit norm. From
Proposition \ref{4..1} it follows that $\norm{e^\varphi a}
=\sup_{\nu\in V} \nu(e^\varphi a)$. For each $\nu\in V$ the
functional $\nu(\:\cdot\:a)$ is positive. And by Proposition
\ref{2..3} the functional $\lambda_\nu(\varphi) =\ln\nu(e^\varphi
a)$ is convex with respect to $\varphi\in\cal C$. Therefore
$\lambda(\varphi) =\sup_{\nu\in V}\lambda_\nu(\varphi)$ is convex
with respect to $\varphi$ too.
\qed
\proofskip

Any linear functional $\nu$ on a positive algebra $\cal A$ we can
multiply by elements $a\in \cal A$ from the left and from the right.
Namely, by definition put $\nu a(\,\cdot\,) =\nu(a\,\cdot\,)$ and
$a\nu(\,\cdot\,) =\nu(\,\cdot\,a)$. If here $\nu$ and $a$ are positive
then the functionals $\nu a$ and $a\nu$ are also positive. For any
$a\in\cal A$ define its spectral radius $r(a) =\lim_{n\to\infty}
\sqrt[\scriptstyle n]{\norm{a^n}}$.

\make{theorem}{4..3}
If\/ $a$ is a positive element of a positive algebra\/ $\cal A$
then there exists a nonzero positive linear functional\/ $\mu\!:
\cal A \to\mathbb R$ such that\/ $a\mu =r(a)\mu$ and\/ $\norm\mu
=\mu(I)$. Similarly, there exists a positive linear functional\/
$\mu'$ such that\/ $\mu'a =r(a)\mu'$ and\/ $\norm{\mu'} =\mu'(I)$.
\end{theorem}

\proof.
Without loss of generality we can assume that $\cal A$ is
complete in norm and that the cone $\cal A_+$ is closed. Consider
$b=r(a)^{-1}a$. The spectral radius of the $b$ is equal to one.
Hence for any $\theta\in(0,1)$ the series $R_\theta =
\sum_{i=0}^\infty \theta^i b^i$ converges in norm and $R_\theta$
is a positive element of $\cal A$. By positiveness of all the
summands in $R_\theta$ the norm of $R_\theta$ does not decrease
as $\theta$ grows. Let $N=\lim_{\theta \nearrow 1}\|R_\theta\|$.
First assume that $N=\infty$. By Proposition~\ref{4..1} for any
$\theta\in(0,1)$ there exists a positive linear functional $\nu$
on $\cal A$ such that $\norm\nu =1$ and $\nu(R_\theta) =
\|R_\theta\|$. Set $\nu_\theta =\|R_\theta\|^{-1}R_\theta\nu$.
Then $\|\nu_\theta\| \le 1$ and at the same time $\nu_\theta(I)
=\norm{R_\theta}^{-1}\nu(R_\theta) =1$. Hence $\|\nu_\theta\|
=1$. Further, $R_\theta -\theta b R_\theta =I$. Therefore the
norm $\|\nu_\theta -\theta b\nu_\theta\| = \|R_\theta\|^{-1}
\|\nu\|$ tends to zero as $\theta\nearrow 1$. Since the unit ball
in the dual space to $\cal A$ is compact in the *-weak topology,
the family $\nu_\theta$ has a limit point as $\theta\nearrow 1$.
This limit point is the desired positive functional $\mu$.

Now assume that $N<\infty$. Then consider the sequence $P_n
=\sum_{i=0}^{2n} i(2n-i)b^i$. Evidently,
\begin{align}
P_n-bP_n &=
\sum_{i=1}^{2n}\bigl[i(2n-i)-(i-1)(2n-i+1)\bigr]b^i = \notag\\
&=\sum_{i=1}^n(2n-2i+1)b^i -\sum_{i=n+1}^{2n}(2i-2n-1)b^i.
\label{4,,1}
\end{align}
From the positiveness of $\cal A$ it follows that for any natural
$n$ and $m>n$ the inequality $\norm{\sum_{i=n}^mb^i}\le N$ holds.
So the norms of both sums in the right hand side of \eqref{4,,1}
do not exceed $2nN$. Consequently, $\|P_n-bP_n\| \le 4nN$. On the
other hand, $\|b^n\|\ge 1$ by the definition of the spectral
radius. Therefore $\|P_n\| \ge \|n^2b^n\|\ge n^2$. The rest of
the proof can be carried out similar to the case $N=\infty$, one
should only replace the family $R_\theta$ by $P_n$ and the
parameter $\theta \nearrow 1$ by $n\to\infty$.

In fact, after the eigenfunctional $\mu$ is found it becomes
clear that the case $N<\infty$ is impossible. Indeed, by
construction $\mu(b^i) =(b^i\mu)(I) =1$ and hence the norm
$\norm{R_\theta}\ge \mu(R_\theta) =\sum_{i=0}^\infty \theta^i$
goes to infinity as $\theta\nearrow 1$.

The existence of the functional $\mu'$ can be proved in the same
way.
\qed

\section{Positive flows}\label{5..}
\setcounter{example}{0}

Suppose $\alpha\!:X\to X$ is a continuous mapping of a Hausdorff
compact set. It generates the shift $\alpha^*(f) =f\circ\alpha$
on $C(X)$. Evidently, this shift is positive and maps the unit
function into itself. Now we shall show that the inverse
statement is true as well.

\begin{proposition}[\mdseries\cite{Homo}]\label{5..1}
If a homomorphism\/ $\delta\!:C(X)\to C(X)$, where\/ $X$ is a
Hausdorff compact set, is positive and maps the unit function
into itself, then there exists a continuous mapping\/ $\alpha
\!:X\to X$ such that\/ $\delta(f)\equiv f\circ\alpha$.
\end{proposition}

This statement is not new, but for the sake of completeness of
the presentation we still give the proof. To each point $x\in
X$ assign the set ${\cal I}_x =\{\,f\in C(X)\mid [\delta(f)](x)
=0\,\}$. Obviously, $\cal I_x$ is a maximal ideal in $C(X)$.
Hence there exists a $y\in X$ such that ${\cal I}_x = \{\,f\in
C(X)\mid f(y) =0\,\}$. Set $\alpha(x) =y$. If $f(y) =c$ then
$f-c\in{\cal I}_x$ and so $[\delta(f)](x) -c =[\delta(f-c)](x)
=0$. Therefore $\delta(f)\equiv f\circ\alpha$. This identity
automatically yields the continuity of $\alpha$.
\qed
\proofskip

Let $\cal C$ be a phase algebra. Any positive homomorphism
$\delta\!:\cal C\to \cal C$ that maps the unit function into
itself we shall call the {\em shift}. Without loss of generality
we can assume that $\cal C$ is dense subalgebra in $C(X)$, where
$X$ is the support of $\cal C$. Then by continuity the shift
$\delta$ extends uniquely from $\cal C$ to $C(X)$ and by
Proposition \ref{5..1} it has the form $\delta(f) =f\circ\alpha$,
where $\alpha$ is a certain continuous mapping of $X$ to itself.
Thus every shift $\delta\!:\cal C\to \cal C$ generates
automatically a dynamical system $(X,\alpha)$.

\proofskip
\mbox{\bfseries Definition.}
A (discrete time) {\em positive flow\/} is a quadruple $(\cal
A,\cal C,E,\delta)$ consisting of a positive Banach (i.\,e.,
complete in norm) algebra $\cal A$ with closed cone $\cal A_+$, a
phase subalgebra $\cal C\subset\cal A$, a fixed element $E\in\cal
A_+$, and a shift $\delta\!:\cal C\to\cal C$ satisfying one of
the two following {\em homological identities\/}: $Ef
=\delta(f)E$ (the covariant form) or $fE = E\delta(f)$ (the
contravariant form) for all $f\in\cal C$.

\proofskip
There is no much difference between the covariant and
contravariant forms of the homological identity. Either of them
one can convert into the other having redefined the
multiplication in $\cal A$ by the rule $a*b =ba$. So in the
sequel we shall use mostly the covariant form.

If a positive flow $(\cal A,\cal C,E,\delta)$ is fixed then the
elements of $\cal A$ we shall call {\em operators\/} and the
distinguished positive element $E$ we shall call the {\em
evolution operator}. The support of the phase algebra $\cal C$ we
shall call the {\em phase space\/} of the positive flow and
denote it by $X$. Once and for all let us assume that the
algebra~$\cal C$ is naturally embedded in $C(X)$. The shift
$\delta$ automatically generates a continuous mapping
$\alpha\!:X\to X$ so that $\delta(f) =f\circ\alpha$ for all $f\in
\cal C\subset C(X)$. Therefore the (covariant) homological
identity takes the form $Ef =(f\circ\alpha)E$. Below, if we wish
to emphasize that the shift $\delta\!:\cal C\to \cal C$ is
generated by a continuous mapping $\alpha\!:X\to X$, we shall use
the notation $\alpha^*$ instead of $\delta$.

Let us give several examples of positive flows.

\make{example}{ex5..1}
Suppose $\alpha\!:X\to X$ is a continuous mapping of a
topological space $X$ to itself. Then one can regard the algebra
of bounded linear operators on $C(X)$ as $\cal A$, the $C(X)$ as
$\cal C$, the weighted shift operator $Ef =e^\varphi
f\circ\alpha$, where $\varphi\in C(X)$ is fixed, as $E$, and the
shift homomorphism $\alpha^*(f) =f\circ\alpha$ as $\delta$.
\end{example}

\make{example}{ex5..2}
Suppose $\alpha\!:X\to X$ is a measurable mapping of a measurable
space supplied with an \hbox{$\alpha$-in}\-var\-iant measure $\mu$.
Then one can regard the algebra of bounded linear operators on
$L^p(X,\mu)$, where $1\le p\le \infty$, as $\cal A$, the algebra
$L^\infty(X,\mu)$ as $\cal C$, the weighted shift operator $Ef
=e^\varphi f\circ\alpha$, provided $\varphi\in L^\infty(X,\mu)$
is fixed, as $E$, and the shift homomorphism $\alpha^*(f)
=f\circ\alpha$ as $\delta$.
\end{example}

\make{example}{ex5..3}
Any stochastic process with discrete time and phase space $X$ one
may think of as the trajectory space $\barX
=\{\,(x_1,x_2,x_3,\dotsc\,) \mid x_i\in X\,\}$ supplied with the
cylinder $\sigma$-algebra and the probability distribution
$\mathbf P$. Then one can regard the algebra of continuous linear
operators on $L^1(\barX,\mathbf P)$ as $\cal A$, the
$L^\infty(\barX,\mathbf P)$ as $\cal C$, and the shift
$Ef(x_1,x_2,x_3,\dotsc\,) =\delta(f)(x_1,x_2,x_3,\dotsc\,)
=f(x_2,x_3,\dotsc\,)$ as $E$ and $\delta$. In essence, this
example is a special case of the previous one.
\end{example}

\make{example}{ex5..4}
Often more profound results come out if one consider a stochastic
process with stationary (invariant under shifts of the time axis)
probability distribution~$\mathbf P$ on the trajectory space
$\barX =X^{\bbb N}$. Then one can define $E$ as the conditional
mathematical expectation $Ef(x_1,x_2,\dotsc\,) =\mathbf E(\,f
(x_0,x_1,x_2,\dotsc\,)\mid x_1,x_2,\dotsc\,)$ on the space
$L^1(\barX,\mathbf P)$ and take the shift $\delta$ on $\cal C
=L^\infty(\barX, \mathbf P)$ from the previous example. In this
case the homological identity takes the contravariant form
$E(\delta(f)g) =fEg$.
\end{example}

\make{proposition}{5..2}
Suppose\/ $(\cal A,\cal C,E,\delta)$ is a positive flow with phase
space\/ $X$. Then there exists a unique algebraic homomorphism\/
$i_{\cal C}\!:C(X)\to \cal A$ that is continuous, injective, and
coincides with the identical embedding on\/ $\cal C\subset C(X)$.
\end{proposition}

\proof.
Consider an arbitrary function $f\in\cal C\subset C(X)$. Put
$t=\max\abs f$. Then $0\le t\pm f\le 2t$ and by monotonicity of
the norm on $\cal A$,
\[
2\norm f =\norm{(t+f)-(t-f)}\le \norm{t+f}+\norm{t-f}\le
\norm{2t}+\norm{2t} =4\max\abs f.
\]
Therefore the identical embedding $i_{\cal C}\!:\cal C\to\cal A$
is continuous with respect to the $\sup$ norm on $\cal C$. By the
continuity it uniquely extends onto $C(X)$. Obviously, the
extended homomorphism is positive. Assume that $f\in C(X)$ and
$i_{\cal C}(f) =0$. Since $\cal C$ is dense in $C(X)$, for any
$\eps>0$ there is a function $g\in\cal C$ such that
$-\varepsilon\le g-f\le\varepsilon$. Applying $i_{\cal C}$ to
these inequalities we obtain $-\varepsilon\le g\le\varepsilon$.
Hence $-2\varepsilon\le f\le 2\varepsilon$. Since $\eps>0$ is
arbitrary, $f$ is zero and so $i_{\cal C}$ is injective.
\qed
\proofskip

Now consider a positive flow $(\cal A,\cal C,E,\alpha^*)$ with
phase space $X$. Suppose $i_{\cal C}\!: C(X)\to\cal A$ is the
homomorphism from Proposition \ref{5..2}. In the sequel we will
identify functions $f\in C(X)$ with their images $i_{\cal C}(f)
\in \cal A$ and will presuppose that $C(X)$ is naturally embedded
into $\cal A$. For each $\varphi\in C(X)$ let us define the
operator $E_\varphi =e^\varphi E\in\cal A$. It is easy to see
that for all $f\in C(X)$ we have the homological identity
\make{equation}{5,,1}
E_\varphi f =e^\varphi Ef =(f\circ\alpha)e^\varphi E =
(f\circ \alpha)E_\varphi.
\end{equation}
Introduce the notation $S_n\varphi
=\varphi+\varphi\circ\alpha+\dots +\varphi\circ\alpha^{n-1}$.
Then \eqref{5,,1} yields $E_\varphi^n = e^\varphi Ee^\varphi
E\dotsm e^\varphi E = e^{S_n\varphi}E^n$. For any function
$\varphi\in C(X)$ define the {\em spectral potential\/}
\make{equation}{5,,2}
\lambda(\varphi) =
\lim_{n\to\infty}\frac{1}{n}\ln\bigl\|E_\varphi^n\bigr\| =
\lim_{n\to\infty}\frac{1}{n}\ln\bigl\|e^{S_n\varphi}E^n\bigr\|.
\end{equation}
In other words, ${\lambda(\varphi)}$ is the logarithm of the
spectral radius of $E_\varphi$.

\make{theorem}{5..3}
The spectral potential of any positive flow is monotone
\textup(if\/ $\varphi\le\psi$ then\/ $\lambda(\varphi)\le
\lambda(\psi)$\textup), additive homogeneous\/ \textup(if\/
$t\in\bbb R$ then\/ $\lambda(\varphi+t) =\lambda(\varphi)
+t$\textup), continuous\/ \textup(satisfies the Lipschitz
condition\/ $|\lambda(\varphi) -\lambda(\psi)|\le
\max|\varphi-\psi|$\textup), strongly\/ $\alpha$-invariant\/
\textup(for all\/ $\varphi,\psi\in C(X)$ satisfies\/
$\lambda(\varphi +\psi) =\lambda(\varphi+\psi\circ
\alpha)$\textup), and convex with respect to\/ $\varphi\in C(X)$.
\end{theorem}

\proof.
The monotonicity and additive homogeneity follow immediately from
\eqref{5,,2} and from the positiveness of $\cal A$. They imply
the Lipschitz condition because $\lambda(\varphi) -
\lambda(\psi)\le\lambda (\psi+\max|\varphi-\psi|) -\lambda(\psi)
=\max|\varphi-\psi|$. Note that always $\bigl|S_n(\varphi+\psi)
-S_n(\varphi+\psi\circ\alpha)\bigr| = |\psi
-\psi\circ\alpha^n|\le 2\sup|\psi|$. Therefore \eqref{5,,2}
yields the strong invariance of the spectral potential. Further,
$S_n\varphi$ depends linearly on $\varphi\in C(X)$. By
Proposition~\ref{4..2} the functional
$\ln\bigl\|e^{S_n\varphi}E^n\bigr\|$ is convex with respect to
$S_n\varphi$ and hence with respect to~$\varphi$ as well. Passing
to the limit in~\eqref{5,,2} we see that $\lambda (\varphi)$ is
also convex.
\qed
\proofskip

It follows from the Lipschitz condition for the spectral
potential that it either is finite on all of $C(X)$ or
$\lambda(\varphi)\equiv-\infty$. In the sequel we always suppose
that the first possibility of the two takes place.

In essence Theorem \ref{5..3} states that the spectral potential
is a particular case of the dynamical potential of a dynamical
system $(X,\alpha)$. So the concepts of dual entropy, equilibrium
measures, and the action functional are already defined for it
(see Section \ref{1..}).

\section{The law of large numbers}\label{6..}

We study the positive flow $(\cal A,\cal C,E,\alpha^*)$ and the
corresponding dynamical system $(X,\alpha)$. Any positive
continuous linear functional on the $\cal A$ we shall call a
{\em state\/} of the flow. The restriction of a state to the
subalgebra $C(X)\subset\cal A$ is also a positive linear
functional. By the Riesz theorem the latter is identified with a
Borel measure on $X$. So every state induces a measure on the
phase space $X$. A state $\nu$ is called normalized if $\nu(1)
=1$. In this case the corresponding measure on $X$ is probability.

Assume that there is some probability measure $P$ on $X$. For any
$f\in C(X)$ the sequence $\{f\circ\alpha^n\}_{n=0}^\infty$ may be
considered as a sequence of random variables with respect to $P$.
This sequence is said to satisfy the law of large numbers if
there exists a constant $a$ such that
\begin{equation*}
\lim_{n\to\infty}P\left\{\abs{\frac{f+f\circ\alpha+\dots +
f\circ\alpha^{n-1}}{n} -a} >\eps\right\} =0 \quad
\hbox{for all $\eps>0$.}
\end{equation*}
This is equivalent to the convergence of $\frac{1}{n}S_nf$ to $a$
with respect to the measure $P$. And if $\frac{1}{n}S_nf$
converges to $a$ almost everywhere, it is said to satisfy the
strong law of large numbers. We will consider one more
modification of the law of large numbers. Assume that $\{P_n\}$
is a sequence of probability distributions on $X$. We shall say
that the sequence $\{f\circ\alpha^n\}$ satisfies the law of large
numbers with respect to the distributions $P_n$ if there exists a
constant $a$ such that
\begin{equation*}
\lim_{n\to\infty}P_n\left\{\abs{\frac{f+f\circ\alpha+\dots +
f\circ\alpha^{n-1}}{n} -a} >\eps\right\} =0 \quad
\hbox{for all $\eps>0$.}
\end{equation*}

Our next purpose is to prove the (strong) law of large numbers
with respect to certain states of a positive flow. To be definite
we consider below only covariant positive flows (although, of
course, each of the statements placed below has a covariant
analog).

Recall that any state $\nu$ of the flow $(\cal A, \cal C,
E,\alpha^*)$ we can multiply by operators $A\in\cal A$ from the
right and from the left: by definition, $\nu A(\,\cdot\,)
=\nu(A\,\cdot\,)$ and $A\nu(\,\cdot\,) =\nu(\,\cdot\,A)$. If $A$
is positive then the products $\nu A$ and $A\nu$ are states too.
A {\em right eigen-state\/} of the positive flow corresponding to
a function $\varphi\in C(X)$ is a state $\nu$ satisfying
$\nu(\,\cdot\,E_\varphi) = e^{\lambda(\varphi)} \nu(\,\cdot\,)$.
Similarly, a {\em left\/} eigen-state of the flow is defined by
the condition $\nu(E_\varphi\,\cdot\,) = e^{\lambda(\varphi)}
\nu(\,\cdot\,)$. If both these conditions take place
simultaneously then the state $\nu$ will be referred to as {\em
two-sided\/} eigen-state. In Theorem \ref{4..3} it was proved
that for any $\varphi\in C(X)$ there exist left and right
eigen-states satisfying $\|\nu\|=\nu(1)$. But the existence of
two-sided eigen-states is still an open problem. Unfortunately, I
do not know any conditions that guarantee the existence of
two-sided eigen-states for an arbitrary positive flow.

A state $\nu$ will be called {\em sensitive\/} with respect to
$\varphi \in C(X)$ (or {\em $\varphi$-sensitive}), if
\begin{equation}\label{6,,1}
\lim_{n\to\infty} \frac{1}{n}\ln \nu\bigl(E_\varphi^n\bigr)
=\lambda(\varphi).
\end{equation}
For example, the left and right eigen-states corresponding to
$\varphi$ are $\varphi$-sensitive. If a set $V\subset C(X)$ is
separable, it is easy to produce a state that is sensitive with
respect to all $\varphi\in V$ (it has the form of a series of
one-sided eigen-states corresponding to the elements of a
countable dense subset in $V$).

\make{proposition}{6..1}
Suppose a state\/ $\nu$ of a positive flow is\/
$\varphi$-sensitive. Then for any\/ $\psi\in C(X)$ and\/ $\eps>0$
there exists so large\/ $N$ that for all\/ $n>N$,
\begin{equation}\label{6,,2}
\frac{\nu\bigl(e^{S_n\psi}E_\varphi^n\bigr)}{\textstyle
\nu\bigl(E_\varphi^n\bigr)} \le
e^{n(\lambda(\varphi+\psi)-\lambda(\varphi) +\eps)}.
\end{equation}
\end{proposition}

\proof.
By the homological identity, $e^{S_n\psi}E_\varphi^n
=e^{S_n(\varphi+\psi)} E^n =(E_{\varphi+\psi})^{n}$. Therefore
$\nu\bigl(e^{S_n\psi}E_\varphi^n\bigr) \le \norm\nu
\norm{(E_{\varphi+\psi})^n}$. Since the spectral radius of
$E_{\varphi+\psi}$ is equal to $e^{\lambda(\varphi+\psi)}$,
combining the last estimate with~\eqref{6,,1} we obtain
\eqref{6,,2}.
\qed
\proofskip

Since the spectral potential $\lambda(\varphi)$ is a particular
case of dynamical potential, we can define on $C^*(X)$ the
corresponding dual entropy $S(\mu) = \inf_{\varphi \in C(X)}
(\lambda(\varphi) -\mu(\varphi))$ and the action functional
$\tau_\varphi(\mu) =\mu(\varphi) +S(\mu) -\lambda (\varphi)$.

Denote by $M(X)$ the set of all probability distributions on $X$.
For each point $x\in X$ define the {\em empirical measures\/}
$\delta_{x,n}\in M(X)$:
\begin{equation*}
\delta_{x,n}(f) =\frac{f(x)+f(\alpha(x))+\dots
+f(\alpha^{n-1}(x))}{n}, \qquad f\in C(X).
\end{equation*}

\make{proposition}{6..2}
Suppose a state\/ $\nu$ of a positive flow is\/
\hbox{$\varphi$-sen}\-sitive. Then for any measure\/ $\mu\in
C^*(X)$ and any number\/ $t>\tau_\varphi(\mu)$ there exists a
small neighborhood\/ $O(\mu)$ in the *-weak topology such that for
all\/ $n$ large enough
\begin{equation}\label{6,,3}
P_n\bigl\{\,x\in X\bigm|\delta_{x,n}\in O(\mu)\,\bigr\}<e^{nt},
\qquad\text{where}\quad P_n =\frac{\nu\bigl(\,\cdot\,E_\varphi^n
\bigr)}{\nu(E_\varphi^n)}.
\end{equation}
\end{proposition}

\proof.
Choose an $\eps>0$ so small that $\tau_\varphi(\mu) <t-\eps$. Then
$S(\mu) =\lambda(\varphi) -\mu(\varphi) +\tau_\varphi(\mu) <
\lambda(\varphi) -\mu(\varphi) +t-\eps$. By the definition of
$S(\mu)$ there exists a function $\psi\in C(X)$ such that
$\lambda(\psi) -\mu(\psi)< \lambda(\varphi) -\mu(\varphi)
+t-\eps$ or, equivalently, $\mu(\psi-\varphi) <\lambda(\psi)
-\lambda(\varphi) -t+\eps$. Set $O(\mu) =\{\,\nu\in C^*(X) \mid
\nu(\psi -\varphi) > \lambda(\psi) -\lambda(\varphi)
-t+\eps\,\}$. Obviously, $S_n(\psi -\varphi)(x) =n\delta_{x,n}
(\psi -\varphi)$. So if $\delta_{x,n}\in O(\mu)$ then $S_n(\psi -
\varphi)(x) -n(\lambda(\psi) -\lambda(\varphi) -t+\eps)>0$. Now
by Proposition~\ref{6..1} for all sufficiently large~$n$ we have
\begin{equation*}
P_n\bigl\{x\in X\bigm|\delta_{x,n}\in O(\mu)\bigr\}
\le P_n\bigl(e^{S_n(\psi-\varphi) -n(\lambda(\psi) -
\lambda(\varphi) -t+\eps)}\bigr) <e^{nt}. \qed
\end{equation*}

Recall that any equilibrium measure corresponding to $\varphi\in
C(X)$ is a subgradient of the functional~$\lambda(\varphi)$ at
the point $\varphi$. In Proposition \ref{1..2} it was proved that
each equilibrium measure is probability and
\hbox{$\alpha$-in}\-var\-iant.

\make{proposition}{6..3}
Assume that an open subset\/ $U\subset M(X)$ contains all
equilibrium measures corresponding to\/ $\varphi\in C(X)$ and that
states\/ $P_n$ have the form\/ $P_n =\nu(\,\cdot\,E_\varphi^n)
\big/ \nu(E_\varphi^n)$, where\/ $\nu$ is some\/
$\varphi$-sensitive state. Then there exists an\/ $\varkappa$
so small that for all large enough\/ $n$ the inequality\/ $P_n
\{\,x\in X \mid \delta_{x,n}\notin U\,\} < e^{-n\varkappa}$ holds.
\end{proposition}

\proof.
Consider any measure $\mu\in M(X)\setminus U$. By Propostition
\ref{1..1} it satisfies $\tau_\varphi(\mu)<0$. By
Propostition~\ref{6..2} there exists a small neighborhood
$O(\mu)$ in $M(X)$ such that $P_n\{\,x\in X\mid \delta_{x,n}\in
O(\mu)\} <e^{n\tau_\varphi (\mu)/2}$ as \hbox{$n\to\infty$}. The
set $M(X)\setminus U$ is compact in the *-weak topology. Hence it
can be covered by finitely many neighborhoods
$O(\mu_1),\,\dotsc,\,O(\mu_k)$ of the form indicated above. Then
$P_n\{\,x\in X \mid \delta_{x,n}\notin U\,\} <\sum_{i=1}^k
e^{n\tau_\varphi(\mu_i)/2}$, provided $n$ is large.
\qed

\make{theorem}{6..4}
Suppose\/ $(\cal A,\cal C,E,\alpha^*)$ is a positive flow with
phase space\/ $X$, and for certain functions\/ $\varphi,f\in
C(X)$ the derivative
\make{equation}{6,,4}
\left.\frac{d\lambda(\varphi+tf)}{dt}\right|_{t=0} =a
\end{equation}
exists, and a state\/ $\nu$ is\/ $\varphi$-sensitive. Then the
sequence\/ $\{f\circ \alpha^n\}$ satisfies the law of large
numbers with respect to the measures\/ $P_n
=\nu(\,\cdot\,E_\varphi^n)\big/\nu(E_\varphi^n)$. Namely,
\make{equation}{6,,5}
\lim_{n\to\infty}P_n\bigl\{\,x\in X\bigm|
\abs{\delta_{x,n}(f) -a} >\eps\,\bigr\} =0
\quad\hbox{for all\,\ $\eps>0$}.
\end{equation}
If, in addition, $\nu$ is a right eigen-state for\/ $E_\varphi$
then the following strong law of large numbers is valid:
\make{equation}{6,,6}
\lim_{n\to\infty}\delta_{x,n}(f) = a \quad
\hbox{for\/ $\nu$-almost all\/ $x$}.
\end{equation}
\end{theorem}

\proof.
It follows from \eqref{6,,4} that any equilibrium measure $\mu$
corresponding to $\varphi$ satisfies $\mu(f) =a$. Let us fix
$\eps>0$ and denote by $U_\eps$ the set of all measures $m\in
M(X)$ such that $\abs{m(f) -a}< \eps$. Thereupon define the
sets $X_{\eps,n} = \{\,x\in X\mid \delta_{x,n} \notin U_\eps\}$. By
Proposition \ref{6..3} there exists $\varkappa>0$ such that $P_n
(X_{\eps,n}) <e^{-n\varkappa}$ for all $n$ large enough. This
estimate implies the first part of Theorem \ref{6..4}.

Now suppose that $\nu$ is a right eigen-state. Without loss of
generality we can assume that it is normalized. Then all $P_n$
coincide with $\nu$. Hence, $\nu \bigl(\bigcup_{n>N}
X_{\eps,n}\bigr) <\sum_{n>N}e^{-n\varkappa} =e^{-N\varkappa}
\big/(e^\varkappa-1)$. Consider the set $X_\eps =
\bigcap_{N>0}\bigcup_{n>N} X_{\eps,n}$. Evidently, $\nu(X_\eps)
=0$. Therefore $\nu\bigl(\bigcup_{\eps>0} X_\eps\bigr) =0$. But
if $x\notin \bigcup_{\eps>0} X_\eps$, then\/ $\lim_{n\to\infty}
\delta_{x,n}(f) =a$.
\qed

\make{corol}{6..5}
If the spectral potential\/ $\lambda(\varphi)$ is G\^ateaux
differentiable at a point\/ $\varphi\in C(X)$, then the laws of
large numbers\/ \eqref{6,,5} and\/ \eqref{6,,6} are valid for all
functions\/ $f\in C(X)$.
\end{corol}

In connection with Theorem \ref{6..4} the question of existence
of the derivative \eqref{6,,4} arises. In convex analysis they
give the following answer to it. Let $\lambda(\varphi)$ be an
arbitrary convex functional on $C(X)$. Then its restriction to
any affine straight line in $C(X)$ can be nondifferentiable on at
most a countable set. If $L$ is a finite-{}dimensional subspace
in $C(X)$ then the restriction of $\lambda(\varphi)$ to $L$ will
be differentiable in the G\^ateaux sense almost everywhere (with
respect to Lebesgue's measure). At last, if $L$ is a separable
vector subspace in $C(X)$, then the restriction of
$\lambda(\varphi)$ to $L$ will be G\^ateaux differentiable on a
residual set (a residual set is a countable intersection of some
open dense subsets of $L$).

Another important question is how the probability distributions
$P_n$ in \eqref{6,,5} might look like. We shall present only one
rather abstract example. Suppose $L$ is a separable vector
subspace in $C(X)$ and $\{\varphi_1,\varphi_2,\ldots\,\}$ is a
countable dense subset in $L$. By Theorem \ref{4..3} for each
$\varphi_i$ there exists a right eigen-state $\nu_i$ such that
$\norm{\nu_i} = \nu_i(I) =2^{-i}$. Define the new state $\nu
=\sum_{i=1}^\infty \nu_i$. It is sensitive with respect to all
$\varphi\in L$. It is easily seen that in this case
\begin{equation*}
\nu_i\bigl(fE_\varphi^n\bigr) =\nu_i\Bigl(fe^{S_n(\varphi
-\varphi_i)} E_{\varphi_i}^n\Bigr) =\nu_i\Bigl(fe^{S_n(\varphi
-\varphi_i)}\Bigr) e^{n\lambda(\varphi_i)}
\end{equation*}
and respectively
\make{equation}{6,,7}
P_n(f) =\frac{\nu\bigl(fE_\varphi^n\bigr)}{\nu\bigl(
E_\varphi^n\bigr)} =
\frac{\sum_{i=1}^\infty \nu_i\bigl(fe^{S_n(\varphi -\varphi_i)}
\bigr) e^{n\lambda(\varphi_i)}}{\sum_{i=1}^\infty
\nu_i\bigl(e^{S_n(\varphi -\varphi_i)}\bigr)
e^{n\lambda(\varphi_i)}}.
\end{equation}

Now let us consider the case of $E_\varphi$ having a two-sided
eigen-state.

\make{proposition}{6..6}
If\/ $\nu$ is a two-sided eigen-state corresponding to\/
$\varphi\in C(X)$ then it induces an\/
\hbox{$\alpha$-in}\-var\-iant measure on\/ $X$.
\end{proposition}

\proof.
This follows from the fact that for any $f\in C(X)$ by the
homological identity we have $\nu(f\circ\alpha) =
e^{-\lambda(\varphi)}\nu((f\circ\alpha) E_\varphi) =
e^{-\lambda(\varphi)}\nu(E_\varphi f) = \nu(f)$.
\qed
\proofskip

\make{proposition}{6..7}
Suppose\/ $\nu$ is a normalized two-sided eigen-state
corresponding to\/ $\varphi\in C(X)$. Then\/~$\nu$ induces an
equilibrium measure on\/ $X$. If, in addition, this measure is
unique equilibrium measure corresponding to\/ $\varphi$ then
it is ergodic.
\end{proposition}

\proof.
Consider the sequence of functionals $\lambda_n(\psi) = n^{-1}\ln
\nu(E_\psi^{n})$. It follows from the definition of the
normalized eigen-state that $\lambda_n(\varphi)
=\lambda(\varphi)$. Recall that by the homological identity
$E_\psi^n =e^{S_n\psi}E^n$. Therefore Proposition~\ref{2..3}
yields the convexity of $\lambda_n(\psi)$. The following
calculation shows that the restriction of $\nu$ to $C(X)$
coincides with the derivative $\lambda'_n(\varphi)$:
\begin{equation*}
\left.\frac{d\lambda_n(\varphi +t\psi)}{dt}\right|_{t=0} =
\frac{1}{n}e^{-n\lambda(\varphi)}\sum_{i=1}^n \nu\bigl(
E_\varphi^{n-i}\psi E_\varphi^i\bigr) =\nu(\psi).
\end{equation*}
Hence $\lambda_n(\varphi+\psi) -\lambda(\varphi)\ge \nu(\psi)$
for all $\psi\in C(X)$. From the last inequality and the
definition of $\lambda_n(\psi)$ it follows that
$\lambda(\varphi+\psi) -\lambda (\varphi) \ge \nu(\psi)$. But
this exactly means that $\nu$ determines an equilibrium
measure. If, in addition, this equilibrium measure is unique for
$\varphi$, then its ergodicity follows from Corollary \ref{6..5}.
\qed
\proofskip

Concluding the section we formulate a problem that seems to be of
great significance. It is well known that in classical
probability theory two theorems have the supreme importance,
namely, the law of large numbers and the central limit theorem.
As we have already seen, the restriction of the spectral
potential to any finite-dimensional subspace $L\subset C(X)$ is
differentiable almost everywhere. And this differentiability
implies the law of large numbers by means of Theorem \ref{6..4}.
But furthermore, there is Alexandroff's theorem \cite{Alexandrov}
asserting that any convex functional on a finite-dimensional
vector space is {\em twice\/} differentiable almost everywhere.
In this connection it would be desirable to know if the central
limit theorem is valid there where the spectral potential is
twice differentiable?

\section{Probabilities of large deviations}\label{7..}

We continue study of the stochastic flow with phase space $X$. Suppose
$L$ is a vector subspace in $C(X)$. For any $\varphi\in L$ let us
define on $L^*$ the {\em conditional action functional\/}
\begin{equation}\label{7,,1}
\tau_\varphi^L(u) =\sup\bigl\{\tau_\varphi(\mu)\bigm|
\left.\mu\right|_L =u\bigr\}, \qquad u\in L^*.
\end{equation}
In the case $L=C(X)$ it coincides with $\tau_\varphi(\mu)$.

\make{proposition}{7..1}
The conditional action functional is upper semicontinuous over $L^*$
\textup(relative to the *-weak topology\textup) and satisfies
\begin{equation}\label{7,,2}
\tau_\varphi^L(u) =\inf_{\psi\in L}\bigl\{\lambda(\psi) -\lambda(\varphi)
-u(\psi-\varphi)\bigr\}, \qquad \varphi\in L,\quad u\in L^*.
\end{equation}
\end{proposition}

\proof.
We have the Young inequality $S(\mu)\le \lambda(\psi) -\mu(\psi)$.
It implies that if $\varphi,\psi\in L$ and $\left.\mu\right|_L =u$
then $\tau_\varphi(\mu) =S(\mu) -\mu(\varphi)
-\lambda(\varphi) \le \lambda(\psi) -\lambda(\varphi)
-u(\psi-\varphi)$. Thus $\tau_\varphi^L (u)$ does not exceed the
right hand side of \eqref{7,,2}.

Conversely, let $t=\inf_{\psi\in L}\{\lambda(\psi)
-\lambda(\varphi) -u(\psi-\varphi)\}$. Then $t+u(\psi) \le
\lambda(\psi) -\lambda(\varphi) +u(\varphi)$ for all $\psi\in
L$. By the Hahn--Banach theorem $u$ extends up to a continuous linear
functional $\mu$ on $C(X)$ such that
$t+\mu(\psi) \le \lambda(\psi) -\lambda(\varphi) +\mu(\varphi)$
for all $\psi\in C(X)$. Therefore,
\begin{equation*}
\tau_\varphi^L(u)\ge \tau_\varphi(\mu) =S(\mu)+\mu(\varphi)
-\lambda(\varphi) =\inf_{\psi\in C(X)}\bigl\{\lambda(\psi)
-\mu(\psi) +\mu(\varphi) -\lambda(\varphi)\bigr\} \ge t.
\end{equation*}
Thus we have proved \eqref{7,,2}. It means that the functional
$\tau_\varphi^L(u) +\lambda(\varphi) -u(\varphi)$, which is equal
to $\inf_{\psi\in L} \{\lambda(\psi) -u(\psi)\}$, only in sign
differs from the Legendre transform of the restriction of
$\lambda(\psi)$ to $L$. Hence it is upper semicontinuous and so
is $\tau_\varphi^L(u)$.
\qed

\make{proposition}{7..2}
If\/ $\tau_\varphi^L(u) >-\infty$ then\/ $u$ belongs to the closure\/
\textup(in the norm of\/ $L^*$\textup) of all subgradients of the
restriction\/ $\left.\lambda(\,\cdot\,)\right|_L$.
\end{proposition}

\proof{}
is just the same as one of Proposition \ref{1..3}.
\qed
\proofskip

Denote by $M_\psi^L$ the set of all subgradients of the restriction
$\left.\lambda(\,\cdot\,)\right|_L$ at the point $\psi \in L$.
Obviously, $M_\psi^L$ is nonempty, convex, and closed in the
*-weak topology. From \eqref{7,,2} it follows
\begin{equation}\label{7,,3}
\tau_\varphi^L(u) \le \lambda(\psi) -\lambda(\varphi)
-u(\psi-\varphi), \qquad u\in L^*.
\end{equation}
Like the Young inequality, \eqref{7,,3} turns into equality iff
$u\in M_\psi^L$.

\make{theorem}{7..3}
Suppose we are given a vector subspace\/ $L\subset C(X)$,
a function\/ $\varphi\in L$, and a certain\/ $\varphi$-sensitive
state\/ $\nu$. Then for any\/ $u\in L^*$ and any number\/
$t_1>\tau_\varphi^L(u)$ there exists a small neighborhood\/
$O(u) \subset L^*$ \textup(in the *-weak topology\textup)
such that for all sufficiently large\/ $n$ we have
\begin{equation}\label{7,,4}
P_n\bigl\{x\in X\bigm| \left.\delta_{x,n}\right|_L\in O(u)\bigr\}
<e^{nt_1}, \qquad\text{where}\quad
P_n =\frac{\nu\bigl(\,\cdot\,E_\varphi^n\bigr)}{\nu(E_\varphi^n)}.
\end{equation}

And if the state\/ $\nu$ is simultaneously sensitive with respect
to two functions\/ $\varphi,\psi\in L$ then for any open subst\/
$U\subset L^*$ containing\/ $M_\psi^L$ and any number\/
$t_2< \inf\bigl\{\tau_\varphi^L(u)\bigm| u\in
M_\psi^L\bigr\}$ we have
\begin{equation}\label{7,,5}
P_n\bigl\{x\in X\bigm| \left.\delta_{x,n}\right|_L\in U\bigr\}
>e^{nt_2}, \qquad \text{where}\quad
P_n =\frac{\nu(\,\cdot\,E_\varphi^n)}{\nu(E_\varphi^n)},
\end{equation}
provided\/ $n$ is large enough.
\end{theorem}

\noindent
\hbox{\bfseries Remark.}
By Proposition \ref{7..2}, if a functional $u\in L^*$ lies at a positive
distance from the set of all subgradients of the restriction
$\left.\lambda(\,\cdot\,)\right|_L$, then $\tau_\varphi^L(u)
=-\infty$. In this case \eqref{7,,4} is valid for any~$t_1<0$.

\make{corol}{7..4}
If under the conditions of Theorem \ref{7..3} the set\/ $M_\psi^L$
consists of only one element\/ $u$ then for any\/ $\eps>0$ there
exists a small neighborhood\/ $O(u)\subset L^*$ such that
\begin{equation}\label{7,,6}
e^{(\tau -\eps)} < P_n\bigl\{x\in X\bigm|
\left.\delta_{x,n}\right|_L \in O(u)\bigr\} <e^{(\tau +\eps)}
\qquad\text{as}\quad n\to \infty,
\end{equation}
where\/ $\tau =\tau_\varphi^L(u) =\lambda(\psi) -\lambda(\varphi) -
u(\psi-\varphi)$.
\end{corol}

The inequalities \eqref{7,,6} represent the main result of the section.
In essence, they mean that the number $\tau_\varphi^L(u)$
determines a rough exponential asymptotics for the probability of
observing the empirical measures whose restriction to $L$ is close
to $u$.

Let us illustrate Corollary \ref{7..4} in the case of finite-dimensional
$L$. Assume that functions $f_1,\dots,f_k\in
C(X)$ are linearly independent and $L=\{\,t^1\!f_1+\dots+t^k\!f_k\mid
f^i\in\mathbb R\,\}$. Denote by $u=(u_1,\dots,u_k)$ the coordinates
in $L^*$ dual to $t= (t^1,\dots,t^k)$. Then formula \eqref{7,,2}
takes the form
\begin{equation*}
\tau_{t_0}^L(u) =\inf_{t\in\mathbb R^k}\bigl\{\Lambda(t) -
\Lambda(t_0) -\langle u,t-t_0\rangle\bigr\},
\end{equation*}
where\ \ $\Lambda(t) =\lambda(t^1\!f_1+\dots+t^k\!f_k)$\ \ and\ \
$t_0^1f_1+\dots+t_0^kf_k =\varphi$.

The functional $\Lambda(t)$ is convex and so is G\^ateaux
differentiable almost everywhere. At the points of differentiability
Corollary \ref{7..4} looks as follows: for any $\eps>0$ there exists
a small $\delta>0$ such that for all sufficiently large $n$ we have
\begin{equation*}
e^{n(\tau-\eps)} <P_n\bigl\{x\in X\bigm|
\abs{\delta_{x,n}(f_i) -u_i}<\delta,\ \ i=1,\,\ldots,\,k\bigr\}
<e^{n(\tau+\eps)},
\end{equation*}
where $u=\Lambda'(t)$\ \ and\ \ $\tau =\tau_{t_0}^L(u) = \Lambda(t)
-\Lambda(t_0) -\langle\Lambda'(t),t-t_0\rangle$.

The finite-dimensional case has the additional advantage that the
state $\nu$, which generates the probability distributions $P_n$,
can be chosen the same for all $t$ and $t_0$. The only condition
for this is sensitivity of $\nu$ with respect to all $\varphi \in
L$.

\proofskip
\proof{}\hbox{\em of Theorem \ref{7..3}}.
We prove the estimate \eqref{7,,4} in just the same way as
\eqref{6,,3}. Choose a small positive $\eps$ such that
$\tau_\varphi^L(u) <t_1-\eps$. By \eqref{7,,2} there exists
$\psi\in L$ satisfying $u(\psi-\varphi) >\lambda(\psi)
-\lambda(\varphi) -t_1 +\eps$. Set $O(u) =\{\,v\in L^*\mid
v(\psi-\varphi) > \lambda(\psi) -\lambda(\varphi) -t_1+\eps\,\}$.
Recall that $S_n(\psi-\varphi)(x) =n\delta_{x,n}(\psi-
\varphi)$. Hence the condition $\left.\delta_{x,n}\right|_L \in
O(u)$ yields the inequality $S_n(\psi-\varphi)(x) -n(\lambda(\psi)
-\lambda(\varphi) -t_1 +\eps) >0$. Then, in view of Proposition
\ref{6..1},
\begin{equation*}
P_n\bigl\{x\in X\bigm| \left.\delta_{x,n}\right|_L \in O(u)\bigr\}
\le P_n\bigl(e^{S_n(\psi-\varphi) -n(\lambda(\psi) -\lambda(\varphi)
-t_1 +\eps)}\bigr) < e^{nt_1},
\end{equation*}
provided $n$ is large. Thus the first part of Theorem \ref{7..3}
is proved.

Now we proceed to the second part. Choose $\eps>0$ such that
$t_2+\eps < \inf\bigl\{\tau_\varphi^L(u)\bigm| u\in
M_\psi^L\bigr\}$, and define the set
\begin{equation}\label{7,,7}
U' =\bigl\{u\in U\bigm| u(\psi -\varphi) <
\lambda(\psi) -\lambda(\varphi) -t_2-\eps\bigr\}.
\end{equation}
Obviously, it is open. And, in view of \eqref{7,,3} and the
choice of $\eps$, it contains $M_\psi^L$.

Let $B$ be the unit ball in $L^*$. By Alaoglu's theorem $B$ is compact
in the *-weak topology. Hence the set $B\setminus U'$ is compact too.
By construction, it has no common points with $M_\psi^L$ and so the
inequality~\eqref{7,,3} is strict on $B\setminus U'$. In addition,
the functional $\tau_\varphi^L(u)$ is upper semicontinuous. Therefore
there exists an $\eps>0$ so small that
\begin{equation}\label{7,,8}
\tau_\varphi^L(u) <\lambda(\psi) -\lambda(\varphi) -u(\psi -
\varphi) -3\eps \qquad\text{for all}\quad u\in B\setminus U'.
\end{equation}

Using the already proved first part of Theorem \ref{7..3}, choose
for each $u\in B \setminus U'$ a neighborhood $O(u)\subset L^*$
so small that
\begin{equation}\label{7,,9}
P_n\bigl\{x\in X\bigm| \left.\delta_{x,n}\right|_L\in O(u)
\bigr\} < e^{n(\tau_\varphi^L(u) +\eps)}
\qquad\text{as}\quad n\to \infty
\end{equation}
and at the same time
\begin{equation}\label{7,,10}
\sup\bigl\{v(\psi -\varphi)\bigm|v\in O(u)\bigr\}
<u(\psi -\varphi) +\eps.
\end{equation}
Choose from these neighborhoods a finite covering
$O(u_1),\,\dots,\,O(u_k)$ of the set $B\setminus U'$. Then
\begin{align*}
P_n\bigl(e^{S_n(\psi -\varphi)}\bigr) &\le
P_n\{\,x\mid \left.\delta_{x,n}\right|_L\in U'\,\}
\sup\bigl\{e^{S_n(\psi-\varphi)(x)}\bigm|
\left.\delta_{x,n}\right|_L \in U'\bigr\} \\[3pt]
&+\sum_{i=1}^k P_n\bigl\{x\bigm| \left.\delta_{x,n}\right|_L
\in O(u_i)\bigr\}
\sup\bigl\{e^{S_n(\psi-\varphi)(x)}\bigm|
\left.\delta_{x,n}\right|_L\in O(u_i)\bigr\}.
\end{align*}
Substituting here the estimates \eqref{7,,7}--\eqref{7,,10} and using
the equality $S_n(\psi-\varphi)(x) =n\delta_{x,n}(\psi -\varphi)$
we obtain
\begin{align}\notag
P_n\bigl(e^{S_n(\psi-\varphi)}\bigr) &\le
P_n\{\,x\mid \left.\delta_{x,n}\right|_L\in U'\,\}e^{n(\lambda(\psi)
-\lambda(\varphi) -t-\eps)} +
\sum_{i=1}^k e^{n(\tau_\varphi^L(u_i) +\eps)}
e^{n(u_i(\psi -\varphi) +\eps)}\\[3pt] \label{7,,11}
&\le e^{n(\lambda(\psi) -\lambda(\varphi))}
\Bigl(P_n\{\,x\mid \left.\delta_{x,n}\right|_L\in U'\,\}
e^{-n(t_2+\eps)} +e^{-n\eps}\Bigr).
\end{align}
On the other hand,
\begin{align}\notag
\lim_{n\to\infty}\frac{1}{n}\ln P_n\bigl(e^{S_n(\psi-\varphi)}\bigr)
& =\lim_{n\to\infty}\frac{1}{n}\ln\frac{\nu\bigl(
e^{S_n(\psi-\varphi)}E_\varphi^n\bigr)}{\nu(E_\varphi^n)}\\[6pt]
& =\lim_{n\to\infty}\frac{1}{n}\ln\frac{\nu(E_\psi^n)}
{\nu(E_\varphi^n)} =\lambda(\psi) -\lambda(\varphi). \label{7,,12}
\end{align}
Comparison of \eqref{7,,11} and \eqref{7,,12} gives \eqref{7,,5}.
\qed

\section{Positive processes}\label{8..}

A (discrete time) {\em positive process\/} is a triple $(\cal
A,\cal C,E)$ consisting of a positive Banach algebra $\cal A$
with closed cone $\cal A_+$, a phase subalgebra $\cal C\subset
\cal A$, and an element $E\in\cal A_+$ named the evolution
operator. For instance, one can take the algebra of bounded
linear operators on the space $L^p(X,\mu)$ as $\cal A$, the
algebra $L^\infty(X,\mu)$ as $\cal C$, and an arbitrary positive
linear operator on $L^p(X,\mu)$ as $E$. The support of the phase
algebra $\cal C$ we will call the phase space of the positive
process and usually denote by $X$. The positive process differs
from the positive flow only by absence of the shift homomorphism
$\delta\!:\cal C\to \cal C$ and associated with it mapping
$\alpha\!:X\to X$. Nevertheless, we shall see soon that to any
positive process one can assign in a standard way a certain
positive flow.

A {\em trajectory\/} of a positive process $(\cal A,\cal C,E)$ is
an arbitrary sequence $\bar x=(x_1,x_2,x_3,\dotsc)$ of points of
the phase space $X$. The set of all trajectories supplied with
the cylinder topology we will denote by $\barX$. By Tychonoff's
theorem $\barX$ is compact.

For every natural $n$ define the subalgebra $C(X^n) \subset
C(\barX)$ as the set of continuous functions of the form
$F(x_1,\dots,x_n)$, which depend only on the first $n$
coordinates of the point $x\in\barX$. The subalgebras $C(X^n)$ are
naturally embedded into each other: $\mathbb R =C(X^0)\subset
C(X^1)\subset C(X^2)\subset\dotsm\,$. By Stone's theorem the
union $\bigcup_n C(X^n)$ is dense in $C(\barX)$.

Consider the left shift $\sigma\!:(x_1,x_2,x_3,\dotsc)
\longmapsto(x_2,x_3,\dotsc)$ on the $\barX$. It generates the
homomorphism of $C(\barX)$ taking each function $F$ to
$F\circ\sigma$. Obviously, $F\circ\sigma^m(x_1,x_2, \dotsc)
=F(x_{1+m},x_{2+m}, \dotsc)$. And if $F\in C(X^n)$ then
$F\circ\sigma^m\in C(X^{n+m})$.

The following object to be defined is the {\em crossed product\/}
of the algebra $C(\barX)$ and the semigroup $\{\,E^n\!\mid n\ge
0\,\}$ generated by the evolution operator $E$. We will use the
notation $\cal B= C(\barX) \times_\sigma \{E^n\}$ for this
crossed product. By definition, the crossed algebra $\cal B$
consists of formal sums $\sum_{i=0}^n F_iE^i$, where $F_i\in
C(\barX)$, which are multiplied by the rule
\make{equation}{8,,1}
\Biggl(\sum_{i=0}^n F_iE^i\!\Biggr)\Biggl(\sum_{j=0}^m G_jE^j
\!\Biggr) =\sum_{i=0}^n\sum_{j=0}^m\bigl(F_i\cdot (G_j\circ
\sigma^i) \bigr)E^{i+j}.
\end{equation}
The elements of $\cal B$ we shall call the {\em crossed
operators}. For any $F\in C(\barX)$ we shall say that the crossed
operator $FE^n$ has {\em degree\/} $n$. From \eqref{8,,1} it
follows that the degree of the product of two crossed operators
is equal to the sum of their degrees. Evidently, the set of all
crossed operators of zero degree forms a subalgebra in $\cal B$
isomorphic to $C(\barX)$. Denote by $\cal B_+$ the set of crossed
operators of the form $\sum_{i=0}^nF_iE^i$, where all the
coefficients $F_i\in C(\barX)$ are nonnegative. We shall suppose
that $\cal B$ is semiordered by the cone $\cal B_+$. A crossed
operator $\sum_{i=0}^n F_iE^i$ will be called {\em proper\/} if
$F_i\in C(X^i)$ for all $i=1,\,\dots,\,n$. The set of all proper
crossed operators we denote by $\mathcal B_p$. Formula
\eqref{8,,1} shows that $\mathcal B_p$ is a subalgebra
in~$\mathcal B$.

\make{theorem}{8..1}
Suppose\/ $(\cal A,\cal C,E)$ is a positive process with phase
space\/ $X$ and\/ $\cal B$ is the corresponding crossed algebra
with the cone\/ $\cal B_+$ and the subalgebra of proper
operators\/ $\cal B_p\subset\cal B$. Then there exists a unique
algebraic homomorphism\/ $\pi\!:\cal B_p\to\cal A$ that is
positive \textup{(}maps the cone\/ $\cal B_p\cap \cal B_+$
to\/~$\cal A_+$\textup{)} and takes each crossed operator of the
form\/ $f_1(x_1)\dotsm f_n(x_n)E^n$, where\/ $f_i\in C(X)$, to
the ordinary product\/ $f_1Ef_2E\dotsm f_nE\in\cal A$. Moreover,
for each function\/ $F\in C(X^n)$ the estimate\/
$\norm{\pi(FE^n)} \le 2\norm{E^n}\sup F$ holds.
\end{theorem}

\proof.
At the beginning let us assume that such a homomorphism does
exist and prove its uniqueness. Consider an arbitrary function
$F\in C (X^n)$. Set $t=\max\abs F$. Then $0\le t\pm F\le 2t$. By
the positiveness of $\pi$ and the monotonicity of the norm on
$\cal A$,
\begin{align}\nonumber
2\norm{\pi(FE^n)} &=\norm{\pi\bigl((t+F)E^n\bigr) -
\pi\bigl((t-F)E^n\bigr)}\\[3pt]\label{8,,2}
&\le \norm{\pi\bigl((t+F)E^n\bigr)} +\norm{\pi\bigl((t-F)E^n\bigr)}
\le 2\norm{\pi(2tE^n)} =4t\norm{E^n}.
\end{align}
Therefore, $\pi(FE^n)$ depends continuously on $F\in C(X^n)$.
Denote by $\cal C_n$ the algebra of functions of the form
$F=\sum_{i=1}^m f^i_1(x_1) \dotsm f^i_n(x_n)$, where $f^i_j\in
C(X)$. By Stone's theorem $\cal C_n$ is dense in~$C(X^n)$. And by
condition $\pi(FE^n)$ is defined uniquely for all $F\in \cal
C_n$. These arguments imply the required uniqueness of~$\pi$.

Now let us prove the existence of $\pi$. By condition,
$\pi(FE^n)$ is already defined for all $F\in\cal C_n$. Assume
that a function $F\in\cal C_n$ is nonnegative. Then by Theorem
\ref{3..4} for any $\eps>0$ the sum $F+\eps$ can be represented
as $F+\varepsilon =\sum_{i=1}^m f^i_1(x_1)\dotsm f^i_n(x_n)$,
where all the functions $f^i_j\in C(X)$ are nonnegative. Hence
$\pi\bigl((F+\varepsilon) E^n\bigr) =\sum_{i=1}^m f^i_1E\dotsm
f^i_nE\in\cal A_+$ and so $\pi(FE^n)\in \cal A_+$ as well. Now
calculation \eqref{8,,2} becomes valid for all $F\in\cal C_n$. It
shows that the linear mapping $F\mapsto \pi(FE^n)$ is continuous
over~$\cal C_n$. By the continuity it extends to $C(X^n)$
preserving the positiveness. Thus we have constructed the linear
mapping $\pi\!:\cal B_p\to\cal A$. It is easily seen that this
$\pi$ is an algebraic homomorphism. Indeed, if $A =f_1(x_1)
\dotsm f_n(x_n)E^n$ and $B =g_1(x_1)\dotsm g_m(x_m)E^m$ then $AB
= f_1(x_1)\dotsm f_n(x_n)\* g_1(x_{n+1}) \dotsm g_m(x_{n+m})
E^{n+m}$ and so $\pi(AB) = f_1E\dotsm f_nE\* g_1E\dotsm g_mE
=\pi(A)\pi(B)$.
\qed
\proofskip

On the crossed algebra $\cal B =C(\barX)\times_\sigma \{E^n\}$
one may define various monotone norms and seminorms (a seminorm
on $\cal B$ is called monotone if $\norm A\le \norm{A+B}$,
provided $A,B\in \cal B_+$). We shall say that a monotone
seminorm $\norm{\,\cdot\,}$ on $\cal B$ is {\em compatible\/}
with the norm on the initial algebra $\cal A$ if $\norm
A=\norm{\pi(A)}$ for any proper positive operator $A\in \cal
B_p\cap \cal B_+$. It turns out that compatible seminorms do
exist and there is a maximal one among them, which can be defined
in the following way. For a positive crossed operator $A\in\cal
B_+$ put
\make{equation}{8,,3}
\norm A =\inf\bigl\{\norm{\pi(A')}\bigm|\, A'\in\cal B_p\quad
\text{and}\quad A'-A\in\cal B_+\bigr\},
\end{equation}
and for the rest $A\in\cal B\,$ let
\make{equation}{8,,4}
\norm A =\inf\bigl\{\norm{A'}+\norm{A''}\bigm|\, A',A''\in\cal B_+
\quad\text{and}\quad A'-A'' =A\bigr\}.
\end{equation}

\make{theorem}{8..2}
Suppose\/ $(\cal A,\cal C,E)$ is a positive process and\/
$\cal B$ is the corresponding crossed algebra with the cone\/
$\cal B_+$ and the subalgebra of proper operators\/ $\cal
B_p\subset\cal B$. Then formulas\/ \eqref{8,,3}, \eqref{8,,4}
determine a monotone seminorm on\/ $\cal B$ compatible with the
norm on\/ $\cal A$. This seminorm is the maximum of all monotone
compatible seminorms on\/ $\cal B$.
\end{theorem}

\proof.
Consider a pair of positive crossed operators $A,B\in \cal B_+$.
Assume that for some proper operators $A',B'\in\cal B_p$ we have
$A'-A\in \cal B_+$ and $B'-B\in\cal B_+$. Then $A',B'\in\cal B_+$
and $A'B'-AB =(A'-A)B' +A(B'-B)\in\cal B_+$. Therefore
\eqref{8,,3} implies $\norm{AB}\le \norm{\pi(A'B')}\le
\norm{\pi(A')}\norm{\pi(B')}$. Passing here to the infimum over
$A',\,B'$ we obtain $\norm{AB}\le \norm A\norm B$. Similarly (and
even a little more easily) it follows from \eqref{8,,3} that
$\norm A+\norm B\ge \norm{A+B}$ and $\norm{A+B}\ge \norm A$.
Further, if $A\in\cal B_+$ then the infimum in \eqref{8,,4} is
attained at $A'=A$ and $A''=0$. Therefore furmulas \eqref{8,,3}
and \eqref{8,,4} give the same value of $\norm A$ for $A\in\cal
B_+$.

Now consider a pair of arbitrary operators $A,B\in\cal B$.
Suppose that some operators $A',A'',B',B'' \in\cal B_+$
satisfy $A'-A'' =A$ and $B'-B'' =B$. Then
\begin{gather*}
A+B=(A'+B')-(A''+B''),\\[3pt]
AB =(A'-A'')(B'-B'') =(A'B' +A''B'')-(A'B''+A''B'),
\end{gather*}
and by virtue of \eqref{8,,4},
\begin{gather}
\norm{A+B}\le\norm{A'+B'}+\norm{A''+B''}\le \norm{A'}
+\norm{A''}+\norm{B'}+\norm{B''}, \label{8,,5}\\[3pt]
\norm{AB}\le \norm{A'B'+A''B''}+\norm{A'B'' +A''B'}\le
\bigl(\norm{A'}+\norm{A''}\bigr)\bigl(\norm{B'}+\norm{B''}\bigr)
\label{8,,6}
\end{gather}
(here we use the inequalities $\norm{C+D}\le \norm C+\norm D$ and
$\norm{CD}\le \norm C\norm D$, which are already proved for
positive operators). Passing in \eqref{8,,5} and \eqref{8,,6} to
the infimums over $A',\,A'',\,B',\,B''$ we get the triangle
inequality $\norm{A+B}\le \norm A+\norm B$ and the estimate for
product $\norm{AB}\le \norm A\norm B$. Eventually we have proved
that \eqref{8,,3} and \eqref{8,,4} determine a monotone seminorm
on $\cal B$. The compatibility of it is evident because in case
$A\in\cal B_p\cap\cal B_+$ the infimum in \eqref{8,,3} is
attained at $A'=A$.

Let us prove that this seminorm is maximal. Consider any other
compatible monotone seminorm $\abs{\,\cdot\,}$ on $\cal B$. If
$A\in\cal B_+$ then, in view of \eqref{8,,3},
\make{equation}{8,,7}
\abs A\le\inf\bigl\{\abs{A'}\bigm|\ A'\in\cal B_p\quad
\text{and}\quad A'-A\in\cal B_+\bigr\} =\norm A,
\end{equation}
and for an arbitrary $A\in\cal B$, in view of \eqref{8,,4},
\eqref{8,,7},
\begin{equation*}
\abs A\le\inf\bigl\{\abs{A'}+\abs{A''}\bigm|\ A',A''\in\cal B_+
\quad\text{and}\quad A'-A'' =A\bigr\}\le\norm A. \qed
\end{equation*}

\make{proposition}{8..3}
Suppose\/ $(\cal A,\cal C,E)$ is a positive process with phase
space\/ $X$, and the seminorm on the corresponding crossed
algebra\/ $\cal B$ is given by\/ \eqref{8,,3}, \eqref{8,,4}. Then
for any function\/ $F\in C(\barX)$ we have\/ $\norm F \le
2\max\abs F\le 2\norm F$.
\end{proposition}

\proof.
The algebra $C(\barX)$ is naturally embedded in $\mathcal B$
(being the subalgebra of crossed operators of zero degree). By
definition, the intersection $C(\barX) \cap\mathcal B_p$ contains
only constant functions. Therefore for nonnegative $F$ furmula
\eqref{8,,3} gives $\norm F =\max\abs F$. For an arbitrary $F\in
C(\barX)$ let $F^+ =\max\{0,F\}$ and $F^- =\max\{0,-F\}$.
Obviously, $F=F^+-F^-$. Hence \eqref{8,,4} gives $\norm F
=\norm{F^+} +\norm{F^-} \le 2\max\abs F$. On the other hand,
$\max\abs F \le \norm{F^+} +\norm{F^-} =\norm F$.
\qed
\proofskip

In the sequel we always suppose that the crossed algebra $\cal B
=C(\barX)\times_\sigma \{E^n\}$ corresponding to the positive
process $(\cal A,\cal C,E)$ is equipped with the maximal
compatible monotone seminorm \eqref{8,,3}, \eqref{8,,4}. In view
of Proposition \ref{8..3} the natural embedding $C(\barX)
\hookrightarrow \mathcal B$ is continuous. Let us form the
quotient of $\cal B$ by the crossed operators with zero seminorm
(i.\,e., all the operators with zero seminorm we regard as
zeros). On this quotient algebra the maximal compatible seminorm
becomes a norm. Complete the quotient algebra in the norm. The
resulting Banach algebra we shall call, as before, the crossed
product and denote it by the same letter $\cal B$. Obviously,
this new $\cal B$ is semiordered (by the closure of the cone
$\cal B_+$) and positive and includes the all of $C(\barX)$.

Formula \eqref{8,,1} yields the homological identity $EF
=(1E)(FE^0) =(F\circ\sigma)E$ for all $F\in C(\barX)$. Denote by
$\sigma^*$ the homomorphism $\sigma^*(F) =F\circ\sigma$, where
$F\in C(\barX)$. Then it follows from the above that the
quadruple $\bigl(\cal B,C(\barX),E,\sigma^* \bigr)$ forms a
positive flow associated with the dynamical system
$(\barX,\sigma)$. We shall call it the {\em suspension\/} over
the positive process $(\cal A,\cal C,E)$. For the suspension
$\bigl(\cal B,C(\barX),E,\sigma^*\bigr)$, just as for any
positive flow, we can define the spectral potential, the
equilibrium measures on $X$, the eigen- and sensitive states. For
them, of course, all the probability-theoretic results proved in
Sections \ref{6..},\,\ref{7..} are valid.

In conclusion we note two more useful properties of the spectral
potential $\lambda(\varphi)$ of the positive flow $\bigl(\cal B,
C(\barX), E,\sigma^*\bigr)$. Firstly, its restriction to the
subalgebra $C(X)\subset C(\barX)$ coincides with the logarithm of
the spectral radius of the operator $e^\varphi E$ in the initial
algebra $\cal A$. This follows from the compatibility of the
norms on $\cal B$ and on $\cal A$. Secondly, the value of
$\lambda (\varphi)$ determined in \eqref{5,,2} does not at all
depend on the choice of a monotone seminorm on the crossed
algebra $\cal B$, provided only this seminorm is compatible with
the norm on $\cal A$. We leave this statement without proof
because it never will be used below.

\section{Supports of positive flows}\label{9..}

Let $(\cal A,\cal C,E,\alpha^*)$ be a covariant positive flow
with phase space $X$. It can be considered as the positive
process $(\cal A,\cal C,E)$. But the latter, in turn, is
identified with the positive flow $\bigl( \cal B,C(\barX),E,
\sigma^*\bigr)$ (suspension). In the current section we
investigate relations between the spectral potentials,
equilibrium measures, and eigen-states of the initial positive
flow and the suspension over it.

A function $f\in C(X)$ will be called {\em unessential\/} for the
positive flow $(\cal A,\cal C,E,\alpha^*)$ if $fE^n =0$ for a
certain $n\in\bbb N$. Evidently, the set of all unessential
functions forms an ideal in $C(X)$. Denote it by $\cal I$. Note
that by the homological identity $EfE^n =(f\circ\alpha) E^{n+1}$.
Hence the ideal $\cal I$ is $\alpha$-invariant.

Define the {\em support\/} of the positive flow as the
annihilator of $\cal I$. In other words, the support is
$\bigcap_{f\in \cal I}f^{-1}(0)$. The $\alpha$-invariance of
$\cal I$ implies $\alpha$-invariance of the support.

\make{proposition}{9..1}
The value of the spectral potential\/ $\lambda(\varphi)$ depends
only on the restriction of\/ $\varphi$ to the support of the
positive flow.
\end{proposition}

\proof.
A function $\psi\in C(X)$ vanishes on the support if and only if
it belongs to the closure of $\cal I$. Therefore it is sufficient
to prove the identity $\lambda(\varphi +\psi) =\lambda(\varphi)$
for all $\psi\in \cal I$. Assume that $\psi E^m =0$. Consider the
operator $E_{\varphi+\psi}^nE^m -E_\varphi^nE^m =\bigl(e^{S_n
(\varphi+ \psi)} -e^{S_n\varphi}\bigr)E^{n+m}$. The difference
$e^{S_n(\varphi+ \psi)} -e^{S_n\varphi}$ can be represented as
$f_n S_n\psi$, where $f_n\in C(X)$. And from the homological
identity it follows that $S_n\psi E^{n+m} =\sum_{i=0}^{n-1}
E^i\psi E^{n+m-i} =0$. Hence $E_{\varphi+\psi}^nE^m
=E_\varphi^nE^m$ for all~$n$. Put $t=\max\abs\varphi$. Since the
norm on $\cal A$ is monotone, $\norm{E_\varphi^n E^m} e^{-mt} \le
\norm{E_\varphi^{n+m}} \le \norm{E_\varphi^n E^m} e^{mt}$. These
inequalities show that $\lambda(\varphi) =\lim_{n\to\infty}n^{-1}
\ln \norm{E_\varphi^n E^m}$. Absolutely in the same way,
$\lambda(\varphi +\psi) =\lim_{n\to\infty}n^{-1}\ln\bigl\|
E_{\varphi+\psi}^n E^m\bigr\|$. Therefore $\lambda(\varphi)
=\lambda(\varphi+\psi)$.
\qed

\make{proposition}{9..2}
Equilibrium measures of a positive flow are concentrated on its
support.
\end{proposition}

\proof.
Each equilibrium measure $\mu$ is a subgradient of the spectral
potential $\lambda(\varphi)$ at a point $\varphi\in C(X)$. It is
sufficient to show that $\mu(\psi) =0$ for all functions $\psi$
that vanish on the support of the flow. By definition of the
subgradient, $\lambda(\varphi+t\psi) -\lambda(\varphi)\ge
\mu(t\psi)$ for all $t\in \mathbb R$. And by Proposition
\ref{9..1} we have $\lambda(\varphi +t\psi) = \lambda(\varphi)$.
Hence $\mu(\psi) =0$.
\qed
\proofskip

We study the positive flow $(\cal A,\cal C,E,\alpha^*)$ with
phase space $X$. Set $\barX =X^{\mathbb N}$ and define the
algebraic homomorphism $\pi\!: C(\barX)\to C(X)$ that takes each
function $F(x_1,x_2,x_3,\dotsc\,)$ to $\pi F(x) =F(x,\alpha(x),
\alpha^2(x),\dotsc\,)$. Consider the crossed product $\cal B
=C(\barX)\times_{\sigma} \{E^n\}$. Extend $\pi$ up to a mapping
$\pi\!:\cal B\to \cal A$ by means of the equality
$\pi\left(\sum_{i=0}^n F_iE^i\right) = \sum_{i=0}^n \pi(F_i)E^i$.
It follows from the homological identity that this extended $\pi$
is also a homomorphism. Which, of course, is positive. On crossed
operators of the form $f_1Ef_2E\ldots f_nE$, where $f_i\in C(X)$,
our extended $\pi$ coincides with the homomorphism that was
defined in Theorem \ref{8..1} and was denoted there by the same
letter $\pi$. And since linear combinations of operators having
the indicated form are dense in the subalgebra of proper crossed
operators $\cal B_p\subset\cal B$, these two homomorphisms
coincide on~$\cal B_p$. Thus we see that the just defined $\pi$
is an extension (from~$\cal B_p$ to $\cal B$) of the homomorphism
from Theorem~\ref{8..1}.

\make{proposition}{9..3}
If a positive flow generates a dynamical system\/
$(X,\alpha)$ then the set of trajectories of the dynamical system
contains the support of the suspension over the flow.
\end{proposition}

\proof.
Suppose a sequence $\bar a=(a_1,a_2,a_3,\dotsc\,)$ is not a
trajectory of the dynamical system $(X,\alpha)$. It is sufficient
to show that $\bar a$ lies out of the support of the suspension
over the initial flow. And for that it is sufficient to find an
unessential (with respect to the suspension) function $F\in
C(\barX)$ such that $F(\bar a)\ne 0$. Fix a number $n$ such that
$\alpha(a_{n-1}) \ne a_{n}$. Choose neighborhoods $O(a_{n-1})$
and $O(a_{n})$ in $X$ such that $\alpha(O(a_{n-1})) \cap O(a_{n})
=\emptyset$. There exist nonnegative functions $f_i \in C(X)$,\
\ $i=n-1,\,n$, vanishing out of $O(a_i)$ and equal to one at
$a_i$. Evidently, $f_{n-1}\cdot(f_{n}\circ\alpha) \equiv 0$.
Consider the function $F(\bar x) = f_{n-1}(x_{n-1})f_{n}(x_{n})$
on $\barX$. By construction, it takes the unit value at~$\bar a$.
Besides, by the compatibility of norms we have $\norm{FE^n}_{\cal
B} = \norm{\pi(FE^n)}_{\cal A} =\norm{E^{n-1}f_{n-1}\cdot(f_{n}
\circ \alpha)E}_{\cal A} =0$. Thus $F$ is unessential.
\qed
\proofskip

Now let $\lambda(\varphi)$ be the spectral potential of a
positive flow $(\cal A,\cal C,E,\alpha^*)$ with phase space $X$
and $\Lambda(\varphi)$ be the spectral potential of the
correspinding suspension $\bigl(\cal B, C(\barX),
E,\sigma^*\bigr)$. Define the embedding $i_\alpha\!: X\to \barX$
that assigns to each point $x\in X$ its trajectory $i_\alpha(x)
=(x,\alpha(x),\alpha^2(x),\dotsc\,)$. This embedding is
associated with the homomorphism $\pi$ (that is $\pi F = F\circ
i_\alpha$).

\make{proposition}{9..4}
The embedding\/ $i_\alpha$ establishes a one-to-one correspondence
between the support of the initial positive flow and the support
of the suspension over it. In addition,
$\Lambda(\varphi)\equiv \lambda(\varphi\circ i_\alpha)$.
\end{proposition}

\proof.
Recall that $C(X)$ is naturally embedded into $C(\barX)$. For any
$\varphi\in C(\barX)$ the difference $\varphi- \varphi\circ
i_\alpha$ vanishes on the set of trajectories $i_\alpha(X)$ and
hence on the support of the suspension too. By Proposition
\ref{9..1} we have $\Lambda(\varphi) \equiv \Lambda(\varphi\circ
i_\alpha)$. And in view of the compatibility of norms the
spectral radius of the operator $e^{\varphi\circ i_\alpha}E$ in
$\cal B$ is equal to the spectral radius of the same operator in
$\cal A$. Therefore, $\lambda(\varphi\circ i_\alpha)
=\Lambda(\varphi\circ i_\alpha) =\Lambda (\varphi)$.

Assume that a trajectory $i_\alpha(a) =(a,\alpha(a), \alpha^2(a),
\dotsc\,)$ lies in the support of the suspension and that a
function $f\in C(X)$ is unessential with respect to the initial
flow. Then $f^2$ is unessential too. Hence for a certain $n$ the
operator $f^2E^n\in \cal A$ is zero. By the compatibility of
norms the operator $f^2E^n\in \cal B_+$ is zero as well.
Therefore~$f^2$ is unessential with respect to the suspension and
so $f^2(i_\alpha(a)) =f^2(a) =0$. From here it follows that~$a$
belongs to the support of the initial flow.

Conversely, assume that $a\in X$ belongs to the support of the
initial flow. We have to prove that if a function $F\in C(\barX)$
is unessential with respect to the suspension then it vanishes
at~$i_\alpha(a)$. Let $FE^n =0$. Then, of course, $F^2E^n =0$.
Consider the function $G(x_1,\dots,x_n) =\max\{\,F^2(x_1,\dots,
x_n,x_{n+1},\dotsc\,)\mid x_{n+1},x_{n+2},\dotsc \in X\,\}$. It
follows from the definition of the maximal compatible seminorm on
the crossed algebra that $\norm{F^2E^n}_{\cal B}
=\norm{\pi(G)E^n}_{\cal A}$. Therefore $\norm{\pi(G)E^n}_{\cal A}
=0$ and hence the function $\pi(G) =G\circ i_\alpha$ is
unessential with respect to the initial flow. Thus
$G(i_\alpha(a)) =0$ and so much the more $F^2(i_\alpha(a)) =0$.
\qed
\proofskip

Proposition \ref{9..4} shows that the embedding $i_\alpha\!:X \to
\barX$ not only establishes a homeomprphism between the support
of the initial positive flow and the support of the suspension
over it but also establishes a one-to-one correspondence between
the sets of equilibrium measures on the both supports. Moreover,
it is easy to see that the embedding $i_\alpha$ conjugates the
mapping $\alpha\!: X\to X$ with the shift $\sigma\!: \barX\to
\barX$, i.\,e., $i_\alpha\circ\alpha =\sigma\circ i_\alpha$.

Finally, let us look at how eigen-states of the initial positive
flow and ones of the corresponding suspension are connected with
each other. Obviously, each state $\nu$ on $\cal A$ induces the
state $\mu =\nu\circ \pi$ on the crossed algebra~$\cal B$. Assume
that $\nu$ is a right eigen-state corresponding to the function
$\pi\varphi = \varphi\circ i_\alpha$, where $\varphi\in
C(\barX)$. Then $\mu =\nu\circ\pi$ is a right eigen-state
corresponding to $\varphi$ because

\begin{equation*}
\mu(\,\cdot\,e^\varphi E) =\nu(\pi(\,\cdot\,)e^{\pi\varphi}E) =
e^{\lambda(\pi\varphi)}\nu(\pi(\,\cdot\,)) =
e^{\Lambda(\varphi)}\mu(\,\cdot\,).
\end{equation*}
On the other hand, the question whether any right eigen-state on
the crossed algebra $\cal B$ can be represented as $\mu =\nu\circ
\pi$, where $\nu$ is a right eigen-state on $\cal A$, still stays
open. In brief, the difficulty of this question is seen from the
fact that the equality $\mu =\nu\circ\pi$ uniquely determines the
restriction of $\nu$ to the subalgebra $\pi(\cal B) \subset\cal
A$ but it is not at all evident that this restriction is positive
(in sense of the cone~$\cal A_+$) or at least continuous with
respect to the norm on $\cal A$.

\section{Examples}\label{10..}
\setcounter{example}{0}

In conclusion we present several examples illustrating the
notions of positive flow, positive process, spectral potential,
dual entropy, and eigen-state.

\make{example}{10..1}
Suppose $\alpha\!:X\to X$ is a continuous mapping of a Hausdorff
compact set $X$. Denote by $\cal A$ the algebra of all continuous
linear operators on $C(X)$. Take $C(X)$ as a phase algebra $\cal
C$. Let an evolution operator $E$ and a homomorphism
$\delta\!:\cal C\to \cal C$ be given by $Ef =\delta(f)
=f\circ\alpha$. For this positive flow the spectral potential
$\lambda(\varphi)$ is the logarithm of the spectral radius of the
weighted shift operator $E_\varphi =e^\varphi E$ that assigns to
each $f\in C(X)$ the function $e^\varphi f\circ\alpha$. It is
well known (see~\cite{Antonevich}, \cite{Kitover},
\cite{Lebedev}, \cite{Lo}) that in this case the spectral
potential has the form $\lambda(\varphi) =\sup_{\mu\in
M_\alpha(X)}\int_X \varphi\, d\mu$, where $M_\alpha(X)$ is the
set of all $\alpha$-invariant Borel probability measures on $X$.
From the remark after Proposition \ref{1..4} it follows that this
spectral potential is maximal among all dynamical potentials with
fixed value $\lambda(0) =0$.

Consider separately the case of zero weight function $\varphi$.
From the form of the functional $\lambda(\varphi)$ we see that it
is G\^ateaux differentiable at the point $\varphi =0$ if and only
if the dynamical system $(X,\alpha)$ has a unique invariant
probability measure $\mu$. In this case $\lambda(\varphi)\equiv
\int_X\varphi\,d\mu$. Dynamical systems with unique invariant
probability measure $\mu$ are called {\em unique ergodic}. In the
unique ergodic dynamical systems all empirical measures converge
to $\mu$ (in the *-weak topology) uniformly with respect to
initial point (see~\cite{Katok}).

Any probability measure $\mu$ on $X$ generates the state $\nu$ on
$\cal A$ that is given by $\nu(A) =\int_X A(1)\,d\mu$. Obviously,
this state satisfies $(E\nu)(A) =\nu(AE) = \nu(A)$. Hence $\mu$
is a right eigen-state corresponding to $\varphi\equiv 0$. And
if, in addition, the measure $\mu$ is invariant then the $\nu$ is
a two-sided eigen-state because $\nu(EA) =\int_X
A(1)\circ\alpha\,d\mu = \int_X A(1)\,d\mu =\nu(A)$.

Now consider an arbitrary function $\varphi\in C(X)$. In the
proof of Theorem \ref{4..3} there was given a way to construct a
right eigen-state corresponding to $\varphi$. Let us try to apply
that construction in our particular case. Evidently, if a linear
operator $A\!:C(X)\to C(X)$ is positive thn $\norm A =\max A(1)$.
For any $\theta<e^{-\lambda (\varphi)}$ consider the operator
$R_\theta = \sum_{n=0}^\infty \theta^n E_\varphi^n$. Its norm is
equal to $\norm{R_\theta} =\max
R_\theta(1) =\max\sum_{n=0}^\infty \theta^n e^{S_n\varphi}$. Let
$\delta_\theta$ be a unit measure concentrated at one of points
where the function $R_\theta(1)$ attains its maximum. Determine
the state $\nu_\theta$ as $\nu_\theta (A) =\delta_\theta(A(1))$.
It has the unit norm and $\nu_\theta (R_\theta)
=\norm{R_\theta}$. In the proof of Theorem \ref{4..3} the right
eigen-state corresponding to $\varphi$ was defined as a limit
point of the family of states $\norm{R_\theta}^{-1}
\nu_\theta(\,\cdot\, R_\theta)$ as $\theta\nearrow
e^{-\lambda(\varphi)}$. The restriction of the state
$\norm{R_\theta}^{-1} \nu_\theta(\,\cdot\,R_\theta)$ to the
algebra $C(X)$ coincides with $\delta_\theta$. Therefore the
restriction of the right eigen-state to~$C(X)$ is a limit point
of the family of measures $\delta_\theta$ as $\theta \nearrow
e^{-\lambda(\varphi)}$. But then this restriction is itself an
atomic measure. Eventually we have proved that for any function
$\varphi\in C(X)$ there exists a right eigen-state such that its
restriction to~$C(X)$ is an atomic probability measure.

Of course, the law of large numbers with respect to an atomic
measure has little significance. But we still can use formula
\eqref{6,,7}, which gives us a sequence of distributed
probabilities $P_n$ such that the law of large numbers is valid
with respect to these $P_n$. For the sake of simplicity let us
consider the case of separable algebra $C(X)$. Suppose
$\{\varphi_1,\varphi_2, \dotsc\,\}$ is a countable dense subset
in $C(X)$. For each $\varphi_i$ there exists a right eigen-state
$\nu_i$ whose restriction to $C(X)$ is concentrated at one point
$x_i$ and has the mass~$2^{-i}$. Then according to \eqref{6,,7}
the measures $P_n$ have the form
\begin{equation*}
P_n(f) =
\frac{\sum_{i=1}^\infty \nu_i\bigl(fe^{S_n(\varphi -\varphi_i)}
\bigr) e^{n\lambda(\varphi_i)}}{\sum_{i=1}^\infty
\nu_i\bigl(e^{S_n(\varphi -\varphi_i)}\bigr)
e^{n\lambda(\varphi_i)}}
=\frac{\sum_{i=1}^\infty
2^{-i}f(x_i)e^{S_n\varphi(x_i) -S_n\varphi_i(x_i) +n\lambda(\varphi_i)}}
{\sum_{i=1}^\infty 2^{-i}e^{S_n\varphi(x_i) -S_n\varphi_i(x_i) +
n\lambda(\varphi_i)}}.
\end{equation*}
They depend on the weight function $\varphi\in C(X)$. If there
exists the G\^ateaux derivative $\lambda'(\varphi) =\mu$ then by
the law of large numbers
\begin{equation*}
\lim_{n\to\infty}P_n\bigl\{\,x\in X\bigm|
\abs{\delta_{x,n}(f) -\mu(f)} >\eps\,\bigr\} =0
\quad\hbox{for all\,\ $\eps>0$}.
\end{equation*}
\end{example}

\make{example}{10..2}
Let $\alpha\!:X\to X$ be a measurable mapping of a measure space
$(X,\mu)$. Consider the shift $E =\alpha^*$ that takes each
function $f$ on $X$ to $Ef =f\circ\alpha$. Assume that there
exists a constant $\cal C$ so large that $\mu(\alpha^{-1}(G))\le
C\mu(G)$ for any measurable subset $G\subset X$. Then the shift
$E$ is continuous on the space $L^1(X,\mu)$ with $\norm E\le C$.
Denote by $\cal A$ the algebra of bounded linear operators on
$L^1(X,\mu)$ and by $\cal C$ the algebra of essentially bounded
functions $L^\infty(X,\mu)$. Then the quadruple $(\cal A,\cal
C,E,\alpha^*)$ forms a positive flow.

Note that the phase space of this flow (that is the support of
$\cal C =L^\infty(X,\mu)$) is not the $X$ but a certain
compactification of $X$. Finite Borel measures on the phase space
are identified with continuous linear functionals on the phase
algebra $\cal C =L^\infty(X,\mu)$ and these functionals, in turn,
are identified with finitely additive measures on $X$ that are
absolutely continuous with respect to $\mu$. For any finitely
additive measure $m$ absolutely continuous with respect to $\mu$
we can define the dual entropy $S(m) = \inf_{\varphi \in \cal
C}\bigl(\lambda(\varphi) -m(\varphi) \bigr)$, where $\lambda
(\varphi)$ is the logarithm of the spectral radius of the
weighted shift operator $E_\varphi =e^\varphi E$. It turns out
that the restriction of $S(m)$ to the set of equilibrium measures
(subgradients of the spectral potential~$\lambda(\varphi)$) can
be computed in a different way avoiding the use of the spectral
radius. Denote by $\goth M_\alpha$ the set of all $\alpha$-invariant
finitely additive probability measures on $X$ that are absolutely
continuous with respect to~$\mu$. We shall consider all possible
finite measurable partitions $D=\{G_1,\dots, G_k\}$ of the space
$X$. Define the functional
\make{equation}{10,,1}
\tau(m) =\inf_{n,D}\sup_{\norm f=1}\sum_{G\in D} m(G)\ln
\frac{\int_G f\circ\alpha^n\,d\mu}{m(G)}, \qquad m\in\goth M_\alpha,
\end{equation}
where the supremum is taken over all nonnegative functions $f\in
L^1(X,\mu)$ with unit norm. In the papers \cite{ABL1},
\cite{ABL2} there was announced the equality $\lambda(\varphi)
=\max_{m\in\goth M_\alpha} \bigl(m(\varphi) +\tau(m)\bigr)$. It
was called there the ``variational principle'' for the spectral
radius of the weighted shift operator. Since the Legendre
transform is involution, this principle implies that the
functionals $S(m)$ and $\tau(m)$ coincide on $\goth M_\alpha$.

Special cases of Example \ref{10..2} are all stationary discrete
time stochastic processes. To verify this fact it is sufficient
to define $X$ as the space of trajectories of a stochastic
process, define $\mu$ as the probability distribution on the $X$,
and define $\alpha$ as the shift taking each trajectory $x(t)$ to
the trajectory~$x(t+1)$. In this case we can define the state
$\nu_0$ on the algebra $\cal A$ of all bounded linear operators
on $L^1(X,\mu)$ by means of the formula $\nu_0(A) =\int_X
A(1)\,d\mu$. The restriction of $\nu_0$ to the phase algebra
$\cal C =L^\infty (X,\mu)$ coincides with $\mu$. Our stochastic
process is stationary iff the measure $\mu$ is
\hbox{$\alpha$-in}\-var\-iant. But the $\alpha$-invariance of
$\mu$ yields that $\nu_0$ is a two-sided eigen-state for the
shift $Ef =f\circ\alpha$ on~$L^1(X,\mu)$.
\end{example}

\make{example}{10..3}
Suppose $\alpha\!:X\to X$ is a continuous and locally expanding
mapping of a metric compact set~$X$. Denote by $H^\gamma(X)$ the
set of real-valued functions on $X$ satisfying the H\"older
condition with exponent $\gamma>0$. For each function $\varphi\in
H^\gamma(X)$ define the following operator $E_\varphi$ on $C(X)$:
\begin{equation*}
E_\varphi f(x) =\sum_{y\in\alpha^{-1}(x)} e^{\varphi(y)}f(y).
\end{equation*}
Usually they call it the Perron--Frobenius operator or the
transfer operator. It is easy to see that it is positive and
satisfies the contravariant homological identity
$E_\varphi\bigl((f\circ\alpha)g\bigr) = fE_\varphi g$. The
operator $E_\varphi$ has a unique positive eigenfunction
$h_\varphi\in H^\gamma(X)$ and a unique positive eigenfunctional
$\mu_\varphi\! :C(X)\to\mathbb R$, and the corresponding
eigenvalue $e^{\lambda_\varphi}$ coincides with the spectral
radius of $E_\varphi$ (see~\cite{P-F1}, \cite{P-F2}). All these
eigen-elements depend analytically on $\varphi \in H^\gamma(X)$
and satisfy
\begin{equation*}
E_\varphi h_\varphi =e^{\lambda_\varphi}h_\varphi,\quad
\mu_\varphi\circ E_\varphi =e^{\lambda_\varphi}\mu_\varphi,\quad
\mu_\varphi(h_\varphi) =1,\quad \mu_\varphi(1) =1.
\end{equation*}

Let us define $\cal A$ as the algebra of bounded linear operators
on $C(X)$ with the transposed multiplication $A*B =B\circ A$,
take $\cal C =H^\gamma(X)$, and define $\delta$ as the shift
$\delta(f) =f\circ\alpha$. Then for any functions $f\in\cal C$
and $g\in C(X)$ we have $(E_\varphi*f)g =fE_\varphi g =E_\varphi
\bigl(\delta(f)g\bigr) = \bigl(\delta(f)*E_\varphi\bigr)g$. So we
get the covariant homological identity $E_\varphi*f
=\delta(f)*E_\varphi$. Hence the quadruple $(\cal A,\cal
C,E_\varphi,\delta)$ forms a covariant positive flow. The number
$\lambda_\varphi$ determined above is its spectral potential. In
our example the evolution operator $E_\varphi$ has a two-sided
eigen-state $\nu_\varphi$ on the algebra $\cal A$. It is given by
$\nu_\varphi(A) =\mu_\varphi(Ah_\varphi)$. The restriction of
$\nu_\varphi$ to $C(X)$ induces an $\alpha$-invariant ergodic
probability measure on $X$. Any function $f\in C(X)$ may be
considered as a random variable with respect to the measure just
defined. In~\cite{P-F3} it was shown that if $f\in H^\gamma(X)$
then the sequence of the random variables $S_nf
=f+f\circ\alpha+\dots+ f\circ\alpha^{n-1}$ satisfies the central
limit theorem and, moreover, Cram\'er's type asymptotics for the
probabilities of large deviations is valid too.
\end{example}

\make{example}{10..4}
Suppose $\Omega$ is an arbitrary measurable space and $B(\Omega)$
is the set of all bounded measurable functions on $\Omega$
supplied with the supremum norm. Assume that $E$ is a positive
linear operator on $B(\Omega)$ such that $E1\equiv 1$. Then $E$
is automatically continuous. Consider the positive process $(\cal
A,B(\Omega),E)$, where $\cal A$ is the algebra of bounded linear
operators on $B(\Omega)$. Denote by $X$ the support of the phase
algebra~$B(\Omega)$. Since $B(\Omega)$ is complete, it is
naturally isomorphic to $C(X)$. The crossed algebra $\cal B$
corresponding to the positive process $(\cal A,B(\Omega),E)$
consists of formal sums $\sum_{i=0}^n F_iE^i$, where the
coefficients $F_i$ are continuous functions on the trajectory
space $\barX =X^{\bbb N}$.

Consider the linear span of all crossed operators of the form
$FE^n$, where $F\in C(\barX)$ is a product $F(x)= f_1(x_1)\dotsm
f_{m}(x_{m})$,\ \ $f_i\in C(X)$. This linear span is dense in the
crossed algebra $\cal B$. Hence any state on $\cal B$ is uniquely
determined by its values on the operators indicated obove. Let
$\mu$ be a positive linear functional on $B(\Omega)$ with $\mu(1)
=1$. Define the state $\nu$ on $\cal B$ by means of the formula
\make{equation}{10,,2}
\nu(FE^n) =\mu(f_1E(f_2E(f_3\dotsc E(f_{m}E(1))...).
\end{equation}
It is easily seen that $\nu$ is a right eigen-state for the
evolution operator $E$. The restriction of $\nu$ to the
subalgebra $C(\barX) \subset\cal B$ generates a Borel probability
measure $\mathbf P$ on $\barX$. In its turn, the linear
functional $\mu$ on $B(\Omega) \cong C(X)$ is identified with a
Borel probability measure on $X$. Further, the evolution
operator~$E$ is naturally transferred from $B(\Omega)$ to $C(X)$.
For each point $x\in X$ consider the positive linear functional
$P_x(f) =Ef(x)$ on $C(X)$. It is also identified with a
probability measure on $X$. Thus we have constructed the Markov
chain with phase space $X$, transition probabilities $P_x$, and
initial probability distribution $\mu$ on $X$. From \eqref{10,,2}
it follows that the corresponding probability measure on the
space of trajectories of this Markov chain coincides with
$\mathbf P$. If we wish, we can suppose that the phase space of
this Markov chain is not $X$, but $\Omega$, although in this case
all the probability measures will be only finitely additive.
\end{example}

\pagebreak


\begin{thebibliography}{99}

\bibitem{Antonevich}
{A.\,B.\ Antonevich.}\ \
Operatoren mit Translation, die durch das Wirken einer kompakten
Lieschen Gruppe erzeugt wird. (German) {\em Sib.\ Mat.\ Zh.,}
{\bfseries 20} (1979), 467-478.
\bibitem{ABL1}
{A.\,B.\ Antonevich, V.\,I.\ Bakhtin, A.\,V.\ Lebedev.}\ \
Variational principle for the spectral radius of weighted shift
and weighted mathematical expectation operators. (Russian)
{\em Dokl.\ Nats.\ Akad.\ Nauk Belarusi,}
{\bfseries 44}, No.~6 (2000), 7--10.
\bibitem{ABL2}
{A.\,B.\ Antonevich, V.\,I.\ Bakhtin, A.\,V.\ Lebedev.}\ \
Theromdynamics and Spectral Radius. {\em Nonlinear Phenomena
in Complex Systems,} {\bfseries 4}, No.~4 (2001), 318--321.
\bibitem{P-F3}
{V.\,I.\ Bakhtin.}\ \
Perron--Frobenius dynamical systems. (Russian)
{\em Dokl.\ Nats.\ Akad.\ Nauk Belarusi,}
{\bfseries 45}, No.~2 (2001), 8--11.
\bibitem{Bishop}
{E.\ Bishop, R.\,R.\ Phelps.}\ \
The support functionals of a convex set, in:
{Proc.\ Sympos.\ Pure Math. 7.} Amer.\ Math.\ Soc., 1963. P.~27--35.
\bibitem{Riesz}
{N.\ Dunford, J.\,T.\ Schwartz.}\ \
Linear Operators. Part I: General Theory. Interscience Publishers,
New York, London, 1958.
\bibitem{Gelfand}
{J.\ Dixmier.}\ \
Les $C^*$-algebres et leurs representations. Paris, 1969.
\bibitem{Alexandrov}
{L.\,C.\ Evans, R.\,F.\ Gariepi.}\ \
Measure theory and fine properties of functions.
CRC Press. Boca Raton, 1992.
\bibitem{P-F1}
{F.\ Hofbauer, G.\ Keller.}\ \
Ergodic properties of invariant measures for piecewice
transformations. {\em Math.~Z.,} {\bfseries 180}, No.~1 (1982),
119--140.
\bibitem{Katok}
{A.\ Katok, B.\ Hasselblatt.}\ \
Introduction to the modern theory of dynamical systems.
Cambridge Univ.\ Press, 1995.
\bibitem{Kitover}
{A.\,K.\ Kitover.}\ \
Spectrum of automorphisms with weight and the Kamowitz-Scheinberg
theorem. {\em Funct.\ Anal.\ Appl.,} {\bfseries 13} (1979), 57-58.
\bibitem{Krasn}
{M.\,A.\ Krasnosel'skij, E.\,A.\ Lifshits, A.\,V.\ Sobolev.}\ \
Positive linear systems. Methods of positive operators.
(Russian) Moscow, 1985.
\bibitem{Lebedev}
{A.\,V.\ Lebedev.}\ \
The invertibility of elements in the $C^*$-algebras generated by
dynamical systems. {\em Russ.\ Math.\ Surv.,} {\bfseries 34}, No.4
(1979), 174-175.
\bibitem{Lo}
{S.\,A.\ L\^o.}\ \
The spectrum of a composition operator in $C(K)$. (Russian)
{\em Izv.\ Akad.\ Nauk BSSR, Ser.\ Fiz.-Mat.\ Nauk},
No.~6 (1979), 44--48.
\bibitem{Homo}
{S.\,A.\ L\^o.}\ \
Weighted shift operators on some Banach spaces of functions.
(Russian) Ph.\ D.\ Thesis. Minsk, 1981.
\bibitem{P-F2}
{D.\ Ruelle.}\ \
The thermodynamic formalism for expanding maps.
{\em Comm.\ Math.\ Phys.,} {\bfseries 125}, No.~2 (1989), 239--262.

\end{thebibliography}
\end{document}